\documentclass[12pt]{amsart}
\usepackage{amsmath,amsthm,enumerate,graphics,amsfonts,latexsym,amsopn,verbatim,amscd,amssymb}
\usepackage{color}
\usepackage{tikz}
\usepackage{amsmath,amssymb,graphicx,tikz,amsthm}
\usetikzlibrary{arrows,matrix}

\theoremstyle{plain}
\newtheorem{thm}{Theorem}[section]

\newtheorem{cor}[thm]{Corollary}
\newtheorem{prop}[thm]{Proposition}

\theoremstyle{definition}

\DeclareMathOperator{\spec}{Spec}
\DeclareMathOperator{\lspec}{L-spec}
\DeclareMathOperator{\qspec}{Q-spec}
\DeclareMathOperator{\cnspec}{CN-spec}
\DeclareMathOperator{\cnlspec}{CNL-spec}
\DeclareMathOperator{\cnqspec}{CNSL-spec}
\DeclareMathOperator{\CNRS}{CNRS}
\DeclareMathOperator{\CNL}{CNL}
\DeclareMathOperator{\CNSL}{CNSL}
\DeclareMathOperator{\CN}{CN}
\DeclareMathOperator{\cnrs}{CNRS}
\DeclareMathOperator{\tr}{tr}
\DeclareMathOperator{\Ecn}{E_{CN}}
\begin{document}

	\title[Common neighborhood Laplacian and signless Laplacian spectra]{Common neighborhood Laplacian and signless Laplacian spectra and energies of commuting graphs}

	\author[F. E. Jannat and R. K. Nath ]{Firdous Ee Jannat and Rajat  Kanti Nath$^*$ }
	\address{Firdous Ee Jannat, Department of Mathematical Science, Tezpur  University, Napaam -784028, Sonitpur, Assam, India.}
	
	\email{firdusej@gmail.com}
	\address{Rajat  Kanti Nath, Department of Mathematical Science, Tezpur  University, Napaam -784028, Sonitpur, Assam, India.}
	
	\email{rajatkantinath@yahoo.com}

\thanks{$^*$Corresponding Author}

\begin{abstract} 
	In this paper, we compute common neighbourhood Laplacian spectrum, common neighbourhood signless Laplacian spectrum and their respective energies of commuting graph of some finite non-abelian  groups including some AC-groups, groups whose  central quotients are isomorphic to $Sz(2)$, $\mathbb{Z}_p\times \mathbb{Z}_p$ or $D_{2m}$. Our findings lead us to conclude that these graphs are CNL (CNSL)-integral. Additionally, we characterize the aforementioned groups such that their commuting graphs   are CNL (CNSL)-hyperenergetic.    
\end{abstract}
	
\thanks{ }
\subjclass[2020]{05C25, 05C50, 15A18}
\keywords{Common Neighborhood; Spectrum; Energy; Commuting Graph}

\maketitle

\section{Introduction} \label{S:intro}
Let $G$ be a finite non-abelian group with center $Z(G)$. The commuting graph of $G$, denoted by $\Gamma_G$ is a simple undirected graph whose vertex set is $G\setminus Z(G)$ and two vertices $x$ and $y$ are adjacent if they commute. This graph was originated from the work of  Brauer and Fowler \cite{BF-55} published in  the year 1955. However, after Neumann's  work \cite{B12} on its complement (that is non-commuting graph) in 1976, the commuting graph became popular. 
Various aspect of commuting graphs of finite non-abelian groups can be found in \cite{AMRR-2006,BBHR-2009,  IJ-2007, MP-2013}. Recently, people have become interested on the spectral aspects of $\Gamma_G$.

Let $\mathcal{G}$ be a simple graph with vertex set $V(\mathcal{G})$ and adjacency matrix $A(\mathcal{G})$. The energy of $\mathcal{G}$ is given by  
$E(\mathcal{G}) :=\sum\limits_{\alpha\in \spec({\mathcal{G}})}|\alpha|$,
where $\spec(\mathcal{G})$ denotes the spectrum of $\mathcal{G}$, that is the set of eigenvalues of $A(\mathcal{G})$ with multiplicites.  In 1978, Gutman \cite{Gutman78} introduced  the notion of $E({\mathcal{G}})$, which has been studied extensively by many mathematicians over the years (see \cite{Gutman-Furtula-2017, GUTMAN-2019*} and the references therein).
 Let $D(\mathcal{G})$ be the degree matrix of  $\mathcal{G}$. The matrices $L(\mathcal{G}) := D(\mathcal{G})-A(\mathcal{G})$ and $Q(\mathcal{G}) := D(\mathcal{G})+A(\mathcal{G})$ are known as  Laplacian  and signless Laplacian matrix of $\mathcal{G}$ respectively.
 Laplacian energy  (denoted by $LE(\mathcal{G})$) and signless Laplacian energy energy  (denoted by $LE^+(\mathcal{G})$) of $\mathcal{G}$ are defined as  
 
   $LE(\mathcal{G}):=\sum\limits_{\beta\in \lspec(\mathcal{G})}\left|\beta-\frac{tr(D(\mathcal{G}))}{|V(\mathcal{G})|}\right|$ and $LE^+(\mathcal{G}):=\sum\limits_{\gamma\in \qspec(\mathcal{G})}\left|\gamma-\frac{tr(D(\mathcal{G}))}{|V(\mathcal{G})|}\right|$,
 
\noindent    where $tr(M)$ denotes the trace of a square matrix $M$; $\lspec(\mathcal{G})$ and $\qspec(\mathcal{G})$ are the sets of eigenvalues (with multiplicity) of $L(\mathcal{G})$ and $Q(\mathcal{G})$, also known as Laplacian and signless Laplacian spectrum of $\mathcal{G}$ respectively. In 2006, Gutman  and Zhou \cite{zhou} introduced  the notion of $LE({\mathcal{G}})$; and in 2008, Gutman et al. \cite{GAVBR} introduced  the notion of $LE^{+}({\mathcal{G}})$.  
   A graph $\mathcal{G}$ is called 
\begin{enumerate} 
 \item hyperenergetic if $E(\mathcal{G}) > E(K_{|V(\mathcal{G})|})$ and borderenergetic if $E(\mathcal{G}) = E(K_{|V(\mathcal{G})|})$. 
\item L-hyperenergetic if $LE(\mathcal{G}) > LE(K_{|V(\mathcal{G})|})$ and L-borderenergetic if $LE(\mathcal{G}) = LE(K_{|V(\mathcal{G})|})$.
\item  Q-hyperenergetic if $LE^+(\mathcal{G}) > LE^+(K_{|V(\mathcal{G})|})$ and Q-borderenergetic if $LE^+(\mathcal{G})$ $ = LE^+(K_{|V(\mathcal{G})|})$.
\end{enumerate} 
%
%
  
The concept of hypererenergetic graph was introduced by Walikar et al.    \cite{Walikar} and Gutman  \cite{Gutman1999}, independently in 1999, while L-hyperenergetic and Q-hyperenergetic graphs were introduced in  \cite{FSN-2020}. The concepts of borderenergetic, L-borderenergetic and Q-borderenergetic graphs were introduced by Gong et al. \cite{GLXGF-2015}, Tura \cite{Tura-2017} and Tao et al. \cite{TH-2018} in the years 2015, 2017 and 2018 respectively.



A graph $\mathcal{G}$ is called integral, L-integral and Q-integral respectively if  $\spec(\mathcal{G})$,  $\lspec(\mathcal{G})$ and  $\qspec(\mathcal{G})$  contain only integers.  The notion of integral graph was introduced by Harary and Schwenk \cite{B10} in 1974. 
 After that Grone and Merris \cite{Grone-94} in 1994 and Simic and
Stanic \cite{Simic-08} in 2008 introduced the notions of L-integral and Q-integral graphs respectively. Integral graphs have some interests for designing the network topology of perfect state transfer networks (see \cite{E1} and the references there in).



In \cite{JR, JR-16}, $\spec(\Gamma_G)$ is computed for various families of finite groups and obtained various groups such that $\Gamma_G$ is integral. In \cite{B4}, Dutta and Nath have computed $\lspec(\Gamma_G)$ and  $\qspec(\Gamma_G)$ for various families of finite groups and obtained various groups such that $\Gamma_G$ is L-integral and Q-integral. In \cite{SND-19,B5,B6}, $E(\Gamma_G)$, $LE(\Gamma_G)$ and $LE^+(\Gamma_G)$ have been computed for many families of finite groups and many finite groups have been obtained for which $\Gamma_G$ is hyperenergetic. It was also shown that $E(\Gamma_G) \leq LE(\Gamma_G)$ for all the groups considered in  \cite{B6}.
In \cite{GAVBR}, it was conjectured that for any graph $\mathcal{G}$,
\begin{equation}\label{E-LE conj}
	E(\mathcal{G}) \leq LE(\mathcal{G}).
\end{equation}
However, \eqref{E-LE conj} (which is also known as E-LE conjecture) was disproved in \cite{Stevanovic,Liu-09}. Since then people wanted to characterize all the graphs subject to the  inequality  \eqref{E-LE conj}.

  In 2011, Alwardi et al. \cite{ASG} introduced the concepts of common neighbourhood spectrum and energy (abbreviated as CN-spectrum and CN-energy)  of a graph $\mathcal{G}$. Let $V(\mathcal{G})=\{v_1,v_2,\ldots,v_n\}$.  
    The common neighbourhood matrix of $\mathcal{G}$, denoted by $CN(\mathcal{G})$,  is a matrix of size $n$  whose $(i, j)$-th entry is given by
  \[
  \CN(\mathcal{G})_{i, j} = \begin{cases} 
  	|C(v_i,v_j)|, & \text{ if } i\neq j\\
  	0, & \text{ otherwise, } 
  \end{cases}
  \]
where $C(v_i,v_j) = \{x \in V(\mathcal{G}) : x \ne v_i, v_j \text{ and   adjacent to both }  v_i, v_j\}$.
%
%
%
  The CN-energy of  $\mathcal{G}$ is defined as
  $E_{CN}(\mathcal{G}) :=\sum\limits_{\delta \in \cnspec(\mathcal{G})}|\delta|$, where $\cnspec(\mathcal{G})$ is the set of eigenvalues of $CN(\mathcal{G})$ with multiplicities.
  A graph $\mathcal{G}$ is called CN-integral if $\cnspec(\mathcal{G})$ contains only integers. Also, $\mathcal{G}$ is called CN-hyperenergetic and CN-borderenergetic (see  \cite{ASG}) if $\Ecn(\mathcal{G})>\Ecn(K_{|V(\mathcal{G})|})$ and  $\Ecn(\mathcal{G})= \Ecn(K_{|V(\mathcal{G})|})$ respectively. 
In \cite{B17,NFDS-2021} Fasfous et al. and Nath et al. discussed various aspects of $\cnspec(\Gamma_G)$ and $E_{CN}(\Gamma_G)$ for several classes of finite non-abelian groups.

In   \cite{FR-2021}, Jannat et al. have introduced the notions of common neighborhood Laplacian spectrum (CNL-spectrum), common neighborhood signless Laplacian spectrum (CNSL-spectrum) and energies corresponding to  those spectra viz. common neighborhood Laplacian energy (CNL-energy) and common neighborhood signless Laplacian energy (CNSL-energy).  The common neighborhood Laplacian matrix (CNL-matrix) and the common neighborhood signless Laplacian matrix (CNSL-matrix) of  $\mathcal{G}$, denoted by $\CNL(\mathcal{G})$ and $\CNSL(\mathcal{G})$, respectively, are given by
\[
\CNL(\mathcal{G}):= \CNRS(\mathcal{G}) - \CN(\mathcal{G}) \text{ and } \CNSL(\mathcal{G}):= \CNRS(\mathcal{G}) + \CN(\mathcal{G}),
\] 
where $\CNRS(\mathcal{G})$ is  a matrix of size $|V(\mathcal{G})|$
whose $(i, j)$-th entry is  given by
 \[
 \CNRS(\mathcal{G})_{i, j} =  \begin{cases}
 	\underset{k = 1}{\overset{|V(\mathcal{G})|}{\sum}}\CN(\mathcal{G})_{i, k}, &   \text{ if }  i = j  \text{ and } i = 1, 2, \dots, |V(\mathcal{G})|\\
 	0, & \text{ if } i \ne j,  
 \end{cases}
 \]

\noindent The set of eigenvalues of  $\CNL(\mathcal{G})$ and $\CNSL(\mathcal{G})$  with multiplicities are called CN-Laplacian spectrum (denoted by $\cnlspec(\mathcal{G})$) and  CN-signless Laplacian spectrum (denoted by $\cnlspec(\mathcal{G})$) of $\mathcal{G}$, respectively. 
Let $\cnlspec(\mathcal{G}) = \{\alpha_1^{a_1}, \,\alpha_2^{a_2}, \dots,\,\alpha_k^{a_k}\}$ and $\cnqspec(\mathcal{G})= \{\beta_1^{b_1},\,\beta_2^{b_2}, \dots,\,\beta_{\ell}^{b_{\ell}}\}$, where $\alpha_1,\,\alpha_2, \dots, \,\alpha_k$ are the distinct eigenvalues of $\CNL(\mathcal{G})$ with corresponding multiplicities  $a_1,\,a_2, \dots,\,a_k$ and $\beta_1,\,\beta_2, \dots,\,\beta_{\ell}$ are the distinct eigenvalues of $\CNSL(\mathcal{G})$ with corresponding multiplicities  $b_1,\,b_2, \dots,\,b_{\ell}$. A graph $\mathcal{G}$ is called CNL (CNSL)-integral if CNL (CNSL)-spectrum  contains only integers.  The  CNL-energy  and  CNSL-energy of $\mathcal{G}$, denoted by $LE_{CN}(\mathcal{G})$ and $LE^+_{CN}(\mathcal{G})$ respectively, are defined as 
	\begin{equation}\label{LEcn}
		LE_{CN}(\mathcal{G}) := \sum_{i=1}^{k} a_i\left|\alpha_i - \Delta_{\mathcal{G}} \right|
	\end{equation}
and
	\begin{equation}\label{LE+cn}
		LE^+_{CN}(\mathcal{G}) := \sum_{i=1}^{\ell} b_i\left| \beta_i - \Delta_{\mathcal{G}} \right|,
	\end{equation}
where $\Delta_{\mathcal{G}}=\frac{\tr(\CNRS(\mathcal{G}))}{|V(\mathcal{G})|}$. 
It was shown, in \cite{FR-2021}, that  
\begin{equation}\label{LEcn-Kn}
	LE_{CN}(K_n) = 	LE^+_{CN}(K_n) = 2(n -1)(n - 2).
\end{equation}
A graph $\mathcal{G}$ is called CNL-hyperenergetic and CNSL-hyperenergetic if $LE_{CN}(\mathcal{G}) > LE_{CN}(K_{|V(\mathcal{G})|}) \text{ and } LE^+_{CN}(\mathcal{G}) > LE^+_{CN}(K_{|V(\mathcal{G})|})$ respectively.
Various aspects of CNL-spectrum, CNL-energy, CNSL-spectrum and CNSL-energy of  graphs, including their relations with other well-known graph energies and Zagreb indices, were discussed in \cite{FR-2021}. 

In this paper, we compute CNL-spectrum, CNSL-spectrum, CNL-energy and CNSL-energy  of commuting graphs of several classes of finite  AC-groups including $QD_{2^n}$ (quasi dihedral group), $PSL(2, 2^k)$ (projective special linear group), $GL(2, q)$ (general linear group where $q >2$ is a prime power), $A(n, v)$, $A(n,p)$ (Hanaki groups), $D_{2m}$ (dihedral group), groups whose central quotient is isomorphic to $Sz(2)$ or $\mathbb{Z}_p\times \mathbb{Z}_p$ or $D_{2m}$ along with some other groups. It follows that the commuting graphs of all the groups considered in our paper are CNL-integral and CNSL-integral. Additionally, we  determine when these graphs are CNL-hyperenergetic and CNSL-hyperenergetic.  
Recall that a group $G$ is called an AC-group if $C_G(x) := \{ y \in G : xy = yx\}$ is abelian for all $x \in G \setminus Z(G)$.

\section{CNL (CNSL)-spectrum and  CNL (CNSL)-energy}
	In this section we have computed the CNL (CNSL)-spectrum and   CNL (CNSL)-energy  of commuting graphs of various groups mentioned above. We start this section with the following result that will be needed for our computations.
\begin{thm}\cite{FR-2021}\label{CN-LE-CNSL-LE^+_Kn}
		Let $\mathcal{G} = l_1 K_{m_1} \cup l_2 K_{m_2}\cup l_3 K_{m_3}$, where $l_iK_{m_i}$ denotes the disjoint union of $l_i$ copies of the complete graphs $K_{m_i}$ on ${m_i}$ vertices for $i = 1, 2, 3$. Then 
		
\noindent		$\cnlspec(\mathcal{G})=\left\lbrace 0^{l_1 + l_2 + l_3}, (m_1(m_1 - 2))^{l_1(m_1 - 1)},\right.$ 
		
			\qquad	\qquad	\qquad	\qquad\qquad\qquad	\qquad $\left.(m_2(m_2 - 2))^{l_2(m_2 - 1)},   (m_3(m_3 - 2))^{l_3(m_3 - 1)}\right\rbrace$ and 
			
\noindent $\cnqspec(\mathcal{G})=\left\lbrace (2(m_1 - 1)(m_1 - 2))^{l_1}, ((m_1 - 2)^2)^{l_1(m_1 - 1)}, (2(m_2 - 1)(m_2 - 2))^{l_2},\right.$
		
\qquad	\qquad	\qquad 	\qquad	 $\left. ((m_2 - 2)^2)^{l_2(m_2 - 1)}, (2(m_3 - 1)(m_3 - 2))^{l_3}, ((m_3 - 2)^2)^{l_3(m_3 - 1)}\right\rbrace$.
\end{thm} 
\subsection{Some families of AC-group}
Here we compute the CNL (CNSL)-spectrum and CNL (CNSL)-energy of  commuting graphs of the quasihedral groups $QD_{2^n}=\{a, b : a^{2^{n - 1}} = b^2 = 1, bab^{-1}=a^{2^{n - 2}}\}$ for  $n\geq 4$,  projective special linear groups $PSL(2, 2^k)$ for $k \geq 2$, general linear groups $GL(2, q)$  for any prime power $q >2$ and the Hanaki groups $A(n, v)$ and $A(n,p)$.
\begin{prop}\label{Quasihedral-theorem}
	The CNL-spectrum, CNSL-spectrum, CNL-energy and CNSL-energy of commuting graphs of the quasihedral groups $QD_{2^n}=\{a, b : a^{2^{n - 1}} = b^2 = 1, bab^{-1}=a^{2^{n - 2}}\}$, where $n\geq 4$, are given by
\[
\cnlspec(\Gamma_{QD_{2^n}})=\{0^{2^{n-1}+1}, ((2^{n - 1} - 2)(2^{n - 1} - 4))^{(2^{n - 1} - 3)}\},
\]
\[
\cnqspec(\Gamma_{QD_{2^n}})=\{0^{2^{n-1}}, (2(2^{n - 1} - 3)(2^{n - 1} - 4))^1 , ((2^{n - 1} - 4)^2)^{(2^{n - 1} - 3)}\},
\]

\centerline{	$LE_{CN}(\Gamma_{QD_{2^n}})=\frac{\left(2^n-8\right) \left(2^n-6\right) \left(2^n-4\right) \left(2^n+2\right)}{8 \left(2^n-2\right)}$ and $LE^+_{CN}(\Gamma_{QD_{2^n}})=\frac{2^{n-3} \left(2^n-8\right) \left(2^n-6\right) \left(2^n-4\right)}{2^n-2}$.}
	
\end{prop}
\begin{proof}
	By \cite[Propositon 2.1]{JR}, we have
	$\Gamma_{QD_{2^n}}=K_{2^{n - 1} - 2}\sqcup 2^{n - 2} K_2.$
	Therefore, by Theorem \ref{CN-LE-CNSL-LE^+_Kn}, we get
	
\centerline{	$\cnlspec(\Gamma_{QD_{2^n}})=\{0^{2^{n - 2} + 1}, ((2^{n - 1} - 2)(2^{n - 1} - 4))^{(2^{n - 1} - 3)}, (2(2-2))^{2^{n - 2}(2-1)} \}$} 

\noindent and 

\vspace{.2cm}

\noindent	$\cnqspec(\Gamma_{QD_{2^n}})=$
	
\noindent	$\{(2(2^{n - 1} - 3)(2^{n - 1} - 4))^1 , ((2^{n - 1} - 4)^2)^{(2^{n - 1} - 3)} , (2(2-1)(2-2))^{(2^{n - 2})} , ((2-2)^2)^{(2^{n - 2})(2-1)}\}$. 

\vspace{.2cm}

\noindent Thus, we get the required	$\cnlspec(\Gamma_{QD_{2^n}})$ and  	$\cnqspec(\Gamma_{QD_{2^n}})$ on simplification.

	\par Here $|V(\Gamma_{QD_{2^n}})|=2^n - 2$ and $\tr(\cnrs(\Gamma_{QD_{2^n}})) = \frac{1}{8} \left(2^n - 8\right) \left(2^n - 6\right) \left(2^n - 4\right)$. Therefore, $\Delta_{\Gamma_{QD_{2^n}}} =  \frac{\left(2^n - 8\right) \left(2^n - 6\right) \left(2^n - 4\right)}{8 \left(2^n - 2\right)}$.  \\
	\noindent\textbf{ CNL-energy calculation:}
	
	We have
	\[ L_1 := |0 - \Delta_{\Gamma_{QD_{2^n}}} | = \bigg| - \frac{\left(2^n - 8\right) \left(2^n - 6\right) \left(2^n - 4\right)}{8 \left(2^n - 2\right)} \bigg|. \]		
Since $n\geq 4$ we get $2^n-8\geq 8$, $2^n-6\geq 10$, $2^n-4\geq 12$ and $2^n-2\geq 14$. So $- \frac{\left(2^n-8\right) \left(2^n-6\right) \left(2^n-4\right)}{8 \left(2^n-2\right)}<0$. Therefore
	 \[L_1=\frac{\left(2^n-8\right) \left(2^n-6\right) \left(2^n-4\right)}{8 \left(2^n-2\right)}.\] 
Also
	 \[L_2 := |\left(2^{n - 1} - 2\right) \left(2^{n - 1} - 4\right) - \Delta_{\Gamma_{QD_{2^n}}} | = \bigg| 2^n(2^{n - 3} - 1) - 1 + \frac{6}{2^n - 2} \bigg|.\] 
Since $n\geq 4$ we have  $2^{n-3}-1\geq 1$ and so  $2^n(2^{n-3}-1)-1 > 0$. Also $2^n-2\geq 14$. Therefore,
\[L_2=2^n(2^{n-3}-1)-1+\frac{6}{2^n-2}.\] 
Hence, by \eqref{LEcn}, we get
\begin{align*}
	LE_{CN}(\Gamma_{QD_{2^n}}) = (2^{n - 1} + 1)L_1 + (2^{n - 1} - 3)L_2 = \frac{\left(2^n - 8\right) \left(2^n - 6\right) \left(2^n - 4\right) \left(2^n + 2\right)}{8 \left(2^n - 2\right)}.
\end{align*}
\noindent \textbf{CNSL-energy calculation:}

We have 
\[ B_1 := | 0 - \Delta_{\Gamma_{QD_{2^n}}} | = \bigg| - \frac{\left(2^n - 8\right) \left(2^n - 6\right) \left(2^n - 4\right)}{8 \left(2^n - 2\right)} \bigg|. \]
Since $n\geq 4$ we get $2^n - 8\geq 8$, $2^n - 6\geq 10$, $2^n - 4\geq 12$ and $2^n - 2\geq 14$. So $ - \frac{\left(2^n - 8\right) \left(2^n - 6\right) \left(2^n - 4\right)}{8 \left(2^n - 2\right)} < 0$. Therefore
\[B_1=\frac{\left(2^n-8\right) \left(2^n-6\right) \left(2^n-4\right)}{8 \left(2^n-2\right)}.\]

We have  
\[ B_2 := \left| 2 \left(2^{n - 1} - 3\right) \left(2^{n - 1} - 4\right) - \Delta_{\Gamma_{QD_{2^n}}} \right| = \left| 15 + 2^n(3\times 2^{n - 3} - 5) + \frac{6}{2^n - 2} \right|. \]
Since $n\geq 4$ we get $3\times 2^{n - 3} - 5 > 0$. Therefore,
\[ B_2 = 15 + 2^n(3\times 2^{n - 3} - 5) + \frac{6}{2^n - 2}. \]
Also,
\[ B_3 := \left| \left(2^{n-1}-4\right)^2 -  \Delta_{\Gamma_{QD_{2^n}}}\right| = \left|7 + 2^{n+1}(2^{n-4}-1) + \frac{6}{2^n-2} \right|. \]
Since $n\geq 4$ we have $2^{n-4}-1>0.$ Therefore, 
\[ B_3 = 7 + 2^{n + 1}(2^{n - 4} - 1) + \frac{6}{2^n - 2}. \]
Hence, by \eqref{LE+cn}, we get
\begin{align*}
	LE^+_{CN}(\Gamma_{QD_{2^ n}}) &= 2^{n - 1}B_1 + 1 \times B_2 + (2^{n - 1} - 3)B_3 = \frac{2^{n - 3} \left(2^n - 8\right) \left(2^n - 6\right) \left(2^n - 4\right)}{2^n - 2}.
\end{align*}
This completes the proof.
\end{proof}

\begin{prop}\label{PSL-theorem}
	The CNL-spectrum, CNSL-spectrum, CNL-energy and CNSL-energy of the commuting graphs of projective special linear groups $PSL(2, 2^k)$, where $k\geq 2$,  are given by

\noindent	$\cnlspec(\Gamma_{PSL(2, 2^k)}) = \left\lbrace (0)^{2^k+2^{2 k}+1}, ((2^k - 1)(2^k - 3))^{(2^k + 1)(2^k - 2)},\right.$
	
	\qquad\qquad\qquad\qquad\qquad\quad $\left. ((2^k - 2)(2^k - 4))^{2^{k - 1} (2^k + 1)(2^k - 3)}, (2^k(2^k - 2))^{2^{k - 1} (2^k - 1)^2}\right\rbrace$,
	
\noindent	$\cnqspec(\Gamma_{PSL(2, 2^k)})$
	 
	 $=\left\lbrace(2 (2^k- 2)(2^k - 3))^{2^k + 1}, ((2^k - 3)^2)^{(2^k + 1)(2^k - 2)}, (2(2^k - 3)(2^k - 4))^{2^{k - 1} (2^k + 1)}, \right.$
	 
	  \quad \,\, $\left. ((2^k - 4)^2)^{2^{k - 1} (2^k + 1)(2^k - 3)}, (2 (2^k - 1)(2^k - 2))^{2^{k - 1} (2^k - 1)},\right. \left. ((2^k - 2)^2)^{2^{k - 1} (2^k - 1)^2}\right\rbrace$,\\	
	$LE_{CN}(\Gamma_{PSL(2, 2^k)}) = \frac{\left(2^k - 2\right) \left(19\times 2^k - 2^{3 k + 2} + 3\times 4^{k + 1} - 7\times 16^k - 5\times 32^k + 3\times 64^k + 6\right)}{8^k - 2^k - 1}$ and
	\begin{align*}
		LE^+_{CN}(\Gamma_{PSL(2, 2^k)}) &= \begin{cases}
												\frac{9260}{59}, &\text{ for } k=2\\
												\frac{- 5\times 2^{k + 2} + 7\times 2^{3 k + 1} + 2^{4 k + 5} - 9\times 4^k + 32^k - 13\times 64^k + 3\times 128^k - 12}{8^k - 2^k - 1}, &\text{ for } k\geq 3.												
											\end{cases}
	\end{align*}

\end{prop}
\begin{proof}
	By \cite[Proposition 2.2]{JR}, we have \[\Gamma_{PSL(2, 2^k)}=(2^k + 1)K_{2^k-1}\sqcup 2^{k - 1} (2^k + 1) K_{2^k - 2} \sqcup 2^{k - 1} (2^k - 1) K_{2^k}.\] 
	Therefore, by Theorem \ref{CN-LE-CNSL-LE^+_Kn}, we get
	
\noindent	$\cnlspec(\Gamma_{PSL(2, 2^k)}) = \left\lbrace (0)^{2^k + 1}, ((2^k - 1)(2^k - 1 - 2))^{(2^k + 1)(2^k - 1 - 1)}, 0^{2^{k - 1} (2^k + 1)},\right.\\ \left. ((2^k - 2)(2^k - 2 - 2))^{2^{k - 1} (2^k + 1)(2^k - 2 - 1)}, 0^{2^{k - 1} (2^k - 1)}, (2^k(2^k - 2))^{2^{k - 1} (2^k - 1)(2^k - 1)}\right\rbrace$ and 
 
\noindent	$\cnqspec(\Gamma_{PSL(2, 2^k)}) = \left\lbrace (2 (2^k-1 - 1)(2^k-1 - 2))^{2^k + 1}, ((2^k-1 - 2)^2)^{(2^k + 1)(2^k-1 - 1)},\right.\\ \left. (2(2^k - 2 - 1)(2^k - 2 - 2))^{2^{k - 1} (2^k + 1)}, ((2^k - 2 - 2)^2)^{2^{k - 1} (2^k + 1)(2^k - 2 - 1)},\right.\\ \left. (2 (2^k - 1)(2^k - 2))^{2^{k - 1} (2^k - 1)}, ((2^k - 2)^2)^{2^{k - 1} (2^k - 1)(2^k - 1)}\right\rbrace$. After simplification, we get the required CNL-spectrum and CNSL-spectrum.
	\par Here $|V(\Gamma_{PSL(2, 2^k)})|=8^k - 2^k - 1$  and  $\tr(\cnrs(\Gamma_{PSL(2, 2^k)}))=\left(2^k - 2\right) \left(5\times 2^k - \right.$ $\left. 3\times 8^k + 16^k + 3\right)$. So, $\Delta_{\Gamma_{PSL(2, 2^k)}} = \frac{\left(2^k - 2\right) \left(5\times 2^k-3\times 8^k + 16^k + 3\right)}{8^k - 2^k - 1}$. \\

	\noindent \textbf{CNL-energy calculation:}
	
	We have
	\[L_1 := \left| 0 - \Delta_{\Gamma_{PSL(2, 2^k)}}\right| = \left| - \frac{\left(2^k - 2\right) \left(5\times 2^k - 3\times 8^k + 16^k + 3\right)}{8^k - 2^k - 1} \right|.\] 
	Since $k\geq 2$ we get $2^k\times 8^k - 3\times 8^k\geq 0$ and $8^k - 2^k - 1 > 0$. So 
		$$- \left(5\times 2^k - 3\times 8^k + 16^k + 3\right)= - \left(2^k\times 8^k - 3\times 8^k + 5\times 2^k + 3\right) < 0.$$ 
		Therefore,
	\[	L_1 = \frac{\left(2^k - 2\right) \left(5\times 2^k - 3\times 8^k + 16^k + 3\right)}{8^k - 2^k - 1}. \]
	We have
	\[L_2 := \left| \left(2^k - 1\right) \left(2^k - 3\right) - \Delta_{\Gamma_{PSL(2, 2^k)}} \right| = \left|\frac{3 + 8 \times 2^k + 4^k (4^k -4 \times 2^k - 2)}{8^k-2^k-1} \right|. \]
Since for all $k \geq 3$, $4^k -4 \times 2^k - 2 > 0$ so $3 + 8 \times 2^k + 4^k (4^k -4 \times 2^k - 2) \geq 0$.  For $k = 2$, $-2^{2 k}+2^{k+3}-2^{3 k+2}-4^k+16^k+3  = 3 > 0$. Again $8^k-2^k-1 > 0$, for all $k\geq 2$. Therefore 
	\[ L_2 = \frac{3 + 8 \times 2^k + 4^k (4^k -4 \times 2^k - 2)}{8^k - 2^k - 1}.\]  

We have  
\[L_3 := \left| \left(2^k - 2\right) \left(2^k - 4\right) - \Delta_{\Gamma_{PSL(2, 2^k)}} \right| = \left| - \frac{ 2 - 5 \times 2^k + 2^k( 2^{2k}(2^k - 1)) }{8^k - 2^k - 1} \right|. \]
Now $2^k-1\geq 3$ and $2^{2k}\geq 16$ so $2 - 5 \times 2^k + 2^k( 2^{2k}(2^k - 1)) > 0$. Therefore $-\frac{ 2 - 5 \times 2^k + 2^k( 2^{2k}(2^k - 1)) }{8^k-2^k-1} < 0.$ Hence
\[L_3=\frac{ 2 - 5 \times 2^k + 2^k( 2^{2k}(2^k - 1)) }{8^k-2^k-1}.\]
Also
\[L_4 := \left| 2^k (2^k - 2) - \Delta_{\Gamma_{PSL(2, 2^k)}} \right| = \left|\frac{ 6 + 9 \times 2^k + 4^k(2 \times 4^k - 4) + 8^k(2 \times 2^k - 7)}{8^k - 2^k - 1} \right|. \] 
Since $2 \times 4^k - 4 > 0$ and $2 \times 2^k - 7 > 0$ we get $6 + 9 \times 2^k + 4^k(2 \times 4^k - 4) + 8^k(2 \times 2^k - 7) > 0$. Therefore
\[L_4=\frac{ 6 + 9 \times 2^k + 4^k(2 \times 4^k - 4) + 8^k(2 \times 2^k - 7)}{8^k - 2^k - 1}.\]
Hence, by \eqref{LEcn}, we get
\begin{align*}
	LE_{CN}&(\Gamma_{PSL(2, 2^k)})= (2^k + 2^{2k} + 1) L_1 + (2^k + 1)(2^k - 2)L_2 + 2^{k-1}(2^k + 1)(2^k - 3)L_3 \\
								&\qquad\qquad\qquad+ 2^{k-1}(2^k-1)^2 L_4\\ 
								&= \frac{\left(2^k - 2\right) \left(19\times 2^k - 2^{3 k + 2} + 3\times 4^{k + 1} - 7\times 16^k - 5\times 32^k + 3\times 64^k + 6\right)}{8^k - 2^k - 1}.
\end{align*}

\noindent \textbf{CNSL-energy calculation:}

We have
\[ B_1 := \left| 2 \left(2^k - 2\right) \left(2^k - 3\right) - \Delta_{\Gamma_{PSL(2, 2^k)}} \right| =  \left| \frac{5\times 2^k + 2^{3 k + 2} + 3\times 4^k - 5\times 16^k + 32^k - 6}{8^k - 2^k - 1} \right|.\]
Now, $\frac{5\times 2^k + 2^{3 k + 2} + 3\times 4^k - 5\times 16^k + 32^k - 6}{8^k - 2^k - 1} > 0$, as $5 \times 2^k - 6 > 0$, for all $k\geq 2$ and $32^k - 5 \times 16^k = 2^k \times 16^k - 5 \times 16^k > 0$, for all $k\geq 3$. For $k=2$, $\frac{5\times 2^k + 2^{3 k + 2} + 3\times 4^k - 5\times 16^k + 32^k - 6}{8^k - 2^k - 1} = \frac{62}{59}>0$. Therefore,  
\[B_1 = \frac{5\times 2^k + 2^{3 k + 2} + 3\times 4^k - 5\times 16^k + 32^k - 6}{8^k - 2^k - 1}. \]
We have
\[B_2 := \left| \left(2^k - 3\right)^2 - \Delta_{\Gamma_{PSL(2, 2^k)}} \right| = \left|\frac{- 3 + 4 \times 2^k - 2^{3k} + 2 \times 2^{3k} - (2^k - 1)2^{3k}}{8^k - 2^k - 1}\right|. \]
As $k\geq 2$, $8^k - 2^k - 1>0$, $2^k-1\geq 3$ and $2^{3k}=2^{2k}\times 2^k\geq 4 \times 2^k$. So $\frac{- 3 + 4 \times 2^k - 2^{3k} + 2 \times 2^{3k} - (2^k - 1)2^{3k}}{8^k - 2^k - 1} < 0$. Therefore
\[B_2 = - \frac{- 3 + 4 \times 2^k - 2^{3k} + 2 \times 2^{3k} - (2^k - 1)2^{3k}}{8^k - 2^k - 1}.\]
We have
{ \small
\[
B_3 := \left| 2 \left(2^k - 3\right) \left(2^k - 4\right) - \Delta_{\Gamma_{PSL(2, 2^k)}} \right|\!=\!\left|\frac{-3\times 2^k+2^{3 k+4}+7\times 4^k-9\times 16^k+32^k\!-\!18}{8^k - 2^k - 1} \right|.
\]}
Now, for $k=2$, $\frac{-3\times 2^k+2^{3 k+4}+7\times 4^k-9\times 16^k+32^k-18}{8^k - 2^k - 1}= - \frac{174}{59} < 0.$ Again $-3\times 2^k+2^{3 k+4}+7\times 4^k-9\times 16^k+32^k-18 = 2^k( 16 \times 2^{2k} - 3) + (7 \times 4^k - 18) + 16(2^k - 9) > 0$, for all $k\geq 4$. For $k=3$, $\frac{-3\times 2^k+2^{3 k+4}+7\times 4^k-9\times 16^k+32^k-18}{8^k - 2^k - 1}= \frac{4502}{503} > 0$. Therefore,
\[B_3=\begin{cases}
		\frac{174}{59}, &\text{ for } k=2\\
	   \frac{-3\times 2^k+2^{3 k+4}+7\times 4^k-9\times 16^k+32^k-18}{8^k - 2^k - 1}, &\text{ for } k\geq 3.
\end{cases}\]  
We have
\[
B_4 := \left| \left(2^k - 4\right)^2 - \Delta_{\Gamma_{PSL(2, 2^k)}} \right| = \left| \frac{-2^k+2^{2 k+1}+9\times 8^k-3\times 16^k-10}{8^k - 2^k - 1} \right|. 
\]
Let $\beta(k):= - 2^k + 2^{2 k + 1} + 9\times 8^k - 3\times 16^k - 10.$ Then  $\beta(k)= - 10 - 2^k + (2^{2k + 1} - 2^{4k}) + (9 \times 8^k - 2\times 2^k \times 8^k) < 0$, as $4k > 2k + 1$ and $2 \times 2^k - 2^{k+1}\geq 9$, for all $k\geq 3$. Again $\beta(2)= - 174 < 0$. Therefore, 
\[B_4 = - \frac{- 2^k + 2^{2 k + 1} + 9\times 8^k - 3\times 16^k - 10}{8^k - 2^k - 1}.\]
We have
{\small
\[
B_5 := \left| 2 \left(2^k - 1\right) \left(2^k - 2\right) - \Delta_{\Gamma_{PSL(2, 2^k)}} \right| = \left| \frac{2 + 9 \times 2^k + 2^{2k}(2^k - 1) 2^k \times 2^k - 4 \times 2^k - 1}{8^k - 2^k - 1} \right|. \]}
Since $k\geq 2$ we get $2^k-1\geq 3$ and $2^k\geq 4$. So $(2^k - 1) 2^k \times 2^k - 4 \times 2^k - 1 > 0.$ Therefore

\[ B_5 = \frac{2 + 9 \times 2^k + 2^{2k}((2^k - 1) 2^k \times 2^k - 4 \times 2^k - 1)}{8^k - 2^k - 1}. \]
Also
\[B_6 := \left| \left(2^k - 2\right)^2 - \Delta_{\Gamma_{PSL(2, 2^k)}} \right| = \left| \frac{2 + 7 \times 2^k + ( ( 2^k - 1) - 3) 8^k + 2^{3k} - 2^{2k + 1}}{8^k - 2^k - 1} \right|. \]
Since $k\geq 2$ we get $(2^k - 1)- 3 = 2^k - 4 \geq 0$ and $2^{3k} \geq 2^{2k + 1}$. Therefore
\[ B_6 = \frac{2 + 7 \times 2^k + ( ( 2^k - 1) - 3) 8^k + 2^{3k} - 2^{2k + 1}}{8^k - 2^k - 1}. \]
Hence, by \eqref{LE+cn}, we get

	\begin{align*}
	LE^+_{CN}(\Gamma_{PSL(2, 2^k)})&= (2^k + 1) B_1 + (2^k + 1)(2^k - 2) B_2 + 2^{k - 1} (2^k + 1) B_3\\
		& + 2^{k - 1} (2^k + 1)(2^k - 3) B_4 + 2^{k - 1} (2^k - 1) B_5 + 2^{k - 1} (2^k - 1)^2 B_6\\
		&= \begin{cases}
				\frac{9260}{59}, &\text{ for } k=2\\
				\frac{- 5\times 2^{k + 2} + 7\times 2^{3 k + 1} + 2^{4 k + 5} - 9\times 4^k + 32^k - 13\times 64^k + 3\times 128^k - 12}{8^k - 2^k - 1}, &\text{ for } k\geq 3.
			\end{cases}
\end{align*}
	Hence the result follows.
\end{proof}

\begin{prop}\label{GL-theorem}
	The CNL-spectrum, CNSL-spectrum, CNL-energy and CNSL-energy of the commuting graphs of general linear groups $GL(2, q)$, where $q=p^n > 2$ and $p$ is a prime, are given by
	
\noindent	$\cnlspec(\Gamma_{GL(2, q)})=\left\lbrace (0)^{q^2 + q + 1}, ((q^2 - 3q + 2)(q^2 - 3q))^{\frac{1}{2}(q^2 + q)(q^2 - 3q + 1)},\right.$
	 
	\qquad \quad $ \left. ((q^2 - q)(q^2 - q - 2))^{\frac{1}{2}(q^2 - q)(q^2 - q - 1)}, ((q^2 - 2q + 1)(q^2 - 2q - 1))^{(q + 1)(q^2 - 2q )}\right\rbrace,$
	
\noindent	$\cnqspec(\Gamma_{GL(2, q)})=\left\lbrace (2(q^2 - 3q + 1)(q^2 - 3q))^{\frac{1}{2}(q^2 + q)}, ((q^2 - 3q)^2)^{\frac{1}{2}(q^2 + q)(q^2 - 3q + 1)},\right.$
	
\qquad \qquad \qquad \qquad  \quad	$\left. (2(q^2 - q - 1)(q^2 - q - 2))^{\frac{1}{2}(q^2 - q)}, ((q^2 - q - 2)^2)^{\frac{1}{2}(q^2 - q)(q^2 - q - 1)},\right.$

\qquad \qquad \qquad \qquad \qquad \quad$ \left. (2(q^2 - 2q)(q^2 - 2q - 1))^{q + 1}, ((q^2 - 2q - 1)^2)^{(q + 1)(q^2 - 2q)}\right\rbrace,$

	$LE_{CN}(\Gamma_{GL(2, q)})=\frac{(q-2) (q-1) q^4 (q+1) (2 q-3) ((q-1) q-1)}{q^3-q-1}$  and \\
	$LE^+_{CN}(\Gamma_{GL(2, q)})=\begin{cases}
		\frac{(q-2) q (q+1) (q (q (q (q (2 q ((q-4) q+6)-7)-7)-1)+11)+2)}{q^3-q-1}, &\text{ if } q\leq 5\\
		\frac{q (q+1) \left((q (q (q (q (2 q-11)+20)-16)+7)+8) q^3-12 q-4\right)}{q^3-q-1}, &\text{ if } q\geq 6.
	\end{cases}$

\end{prop}
\begin{proof}
	By \cite[Proposition 2.3]{JR}, we have 
	\[\Gamma_{GL(2, q)}=\frac{q(q + 1)}{2}K_{q^2 - 3q + 2}\sqcup \frac{q(q - 1)}{2}K_{q^2 - q}\sqcup (q + 1)K_{q^2 - 2q + 1}.\]
	Therefore, by Theorem \ref{CN-LE-CNSL-LE^+_Kn}, we get
	
\noindent	$\cnlspec(\Gamma_{GL(2, q)})=\left\lbrace (0)^{\frac{q(q + 1)}{2}}, ((q^2 - 3q + 2)(q^2 - 3q + 2 - 2))^{(\frac{q(q + 1)}{2})(q^2 - 3q + 2 - 1)},\right.$

\qquad\qquad\qquad\qquad\qquad $\left. (0)^{\frac{q(q - 1)}{2}}, ((q^2 - q)(q^2 - q - 2))^{(\frac{q(q - 1)}{2})(q^2 - q - 1)}, (0)^{q + 1},\right.$ 

\qquad\qquad\qquad\qquad \qquad\qquad \quad$\left. ((q^2 - 2q + 1)(q^2 - 2q + 1 - 2))^{(q + 1)(q^2 - 2q + 1 - 1)}\right\rbrace$ and 
		
\noindent			$\cnqspec(\Gamma_{GL(2, q)})=\left\lbrace (2(q^2 - 3q + 2 - 1)(q^2 - 3q + 2 - 2))^{\frac{q(q + 1)}{2}},\right.$

\qquad\quad\quad$\left. ((q^2 - 3q + 2 - 2)^2)^{(\frac{q(q + 1)}{2})(q^2 - 3q + 2 - 1)}, (2(q^2 - q - 1)(q^2 - q - 2))^{\frac{q(q - 1)}{2}},\right.$

\qquad\quad\quad$\left. ((q^2 - q - 2)^2)^{(\frac{q(q - 1)}{2})(q^2 - q - 1)}, (2(q^2 - 2q + 1-1)(q^2 - 2q + 1-2))^{q + 1},\right.$

\qquad\qquad\qquad\qquad \qquad\qquad\qquad\qquad\qquad\quad\,\,\,$\left. ((q^2 - 2q + 1-2)^2)^{(q + 1)(q^2 - 2q + 1 - 1)}\right\rbrace.$

\noindent	Thus we get $\cnlspec(\Gamma_{GL(2, q)})$ and $\cnqspec(\Gamma_{GL(2, q)})$ on simplification.
	\par Here $|V(\Gamma_{GL(2, q)})|=(q-1) (q^3 - q - 1)$ and $\tr(\cnrs(\Gamma_{GL(2, q)}))=(q-2) (q-1) q (q+1) \left((q-2) (q-1) q^2+1\right)$. So, $\Delta_{\Gamma_{GL(2, q)}} =  \frac{(q - 2) q (q + 1) \left((q - 2) (q - 1) q^2 + 1\right)}{q^3 - q - 1}$. 
	
	\noindent \textbf{CNL-energy calculation:} 
	
	We have
	\[L_1 := \left| 0 - \Delta_{\Gamma_{GL(2, q)}} \right| = \left| - \frac{(q - 2) q (q + 1) ((q - 2) (q - 1) q^2 + 1)}{q^3 - q - 1}\right|. \]
	Note that $- (q - 2) q (q + 1) ((q - 2) (q - 1) q^2 + 1) < 0$ and $q^3 - q - 1 > 0$ since $q > 2$. Therefore
	\[L_1 = \frac{(q - 2) q (q + 1) ((q - 2) (q - 1) q^2 + 1)}{q^3 - q - 1}.\]
	We have
	\[ L_2 := \left| \left(q^2 - 3 q + 2\right) \left(q^2 - 3 q\right) - \Delta_{\Gamma_{GL(2, q)}} \right| = \left| - \frac{(q - 2) q \left((q (2 q - 3) - 1) q^2 + 4\right)}{q^3 - q - 1} \right|.\]
	Note that $ - (q - 2) q \left((q (2 q - 3) - 1) q^2 + 4\right) < 0$ and  $q^3 - q - 1 > 0$ since $q > 2$. Therefore,
	\[L_2 = \frac{(q - 2) q \left((q (2 q - 3) - 1) q^2 + 4\right)}{q^3 - q - 1}.\]
	We have
	\[ L_3 := \left| \left(q^2 - q\right) \left(q^2 - q - 2\right) - \Delta_{\Gamma_{GL(2, q)}} \right| = \left| \frac{(q-2) q^3 (q+1) (2 q-3)}{q^3-q-1} \right|.\]
	Note that $(q-2) q^3 (q+1) (2 q-3) > 0$ since $q > 2$. Therefore, 
	\[L_3 = \frac{(q-2) q^3 (q+1) (2 q-3)}{q^3-q-1}. \]
	Also
	\[ L_4 := \left| \left(q^2 - 2 q + 1\right) \left(q^2 - 2 q - 1\right) - \Delta_{\Gamma_{GL(2, q)}} \right| = \left| \frac{1 - q (-3 + q (3 + (-2 + q) q))}{q^3-q-1} \right|. \]
	Note that $1 - q (-3 + q (3 + (-2 + q) q)) < 0$ since $q > 2$. Therefore,
	\[ L_4 = \frac{q (-3 + q (3 + (-2 + q) q)) - 1}{q^3-q-1}. \]
	Hence, by \eqref{LEcn}, we get
	\begin{align*}
		LE_{CN}(\Gamma_{GL(2, q)}) &= (q^2 + q + 1) L_1 + \frac{q(q + 1)}{2}(q^2 - 3q + 1) L_2 + \frac{q(q - 1)}{2}(q^2 - q - 1) L_3 \\
									&\quad\,+ (q + 1)(q^2 - 2q ) L_4\\
								   &= \frac{(q-2) (q-1) q^4 (q+1) (2 q-3) ((q-1) q-1)}{q^3-q-1}.
	\end{align*}

	\noindent \textbf{CNSL-energy calculation: }
	
	We have
	\[ B_1 := \left| 2 \left(q^2 - 3 q + 1\right) \left(q^2 - 3 q\right) - \Delta_{\Gamma_{GL(2, q)}} \right| = \left| \frac{q \left(q \left((q - 5) (q - 3) q^3 - 5 q - 13\right) + 8\right)}{q^3 - q - 1} \right|. \]
	For $q \geq 6$ we have  
	$q^3 (q - 5)(q - 3) \geq 648$. Therefore,  $q^2 (q - 5)(q - 3) q - 5q - 13 > 0$ and so $q(q ((q - 5) (q - 3) q^3 - 5 q - 13) + 8) > 0.$ Again, for $q \leq 5$ we have $q^3 (q - 5)(q - 3) \leq 0$ and so $q(q ((q - 5) (q - 3) q^3 - 5 q - 13) + 8) < 0$. Therefore, 
	\[B_1 = \begin{cases}
		- \frac{q \left(q \left((q - 5) (q - 3) q^3 - 5 q - 13\right) + 8\right)}{q^3 - q - 1},~~~~&\text{ for } q \leq 5\\
		\frac{q \left(q \left((q - 5) (q - 3) q^3 - 5 q - 13\right) + 8\right)}{q^3 - q - 1},~~~~~~~~&\text{ for } q \geq 6.
	\end{cases}\]
We have
	\[ B_2 := \left| \left(q^2 - 3 q\right)^2 - \Delta_{\Gamma_{GL(2, q)}} \right| = \left| \frac{q \left(- 2 q^5 + 5 q^4 + q^3 - 8 q + 2\right)}{q^3 - q - 1} \right|. \]
	Since $q \geq 3$ we get $q \left(- 2 q^5 + 5 q^4 + q^3 - 8 q + 2\right)= 2 - 8q + q^3(1 + q(5 - 2q)) < 0$ and $q^3 - q - 1 > 0$. Therefore, 
	\[B_2 = - \frac{q \left(- 2 q^5 + 5 q^4 + q^3 - 8 q + 2\right)}{q^3 - q - 1}.\]
	We have
\begin{align*}
	 B_3 :=& \left| 2 \left(q^2 - q - 1\right) \left(q^2 - q - 2\right) - \Delta_{\Gamma_{GL(2, q)}} \right|\\	 
	 &\qquad\qquad\qquad\qquad\qquad = \left| \frac{(q - 2) (q + 1) \left(q^5 + q^4 - 6 q^3 + 3 q + 2\right)}{q^3 - q - 1} \right|. 
 \end{align*}
	Since $q \geq 3$ we get $(q - 2) (q + 1) (q^5 + q^4 - 6 q^3 + 3 q + 2)= (q - 2) (q + 1) (2 + 3q + q^4 + q^3(q^2 - 6))  > 0$ and $q^3 - q - 1 > 0$. Therefore,
	\[B_3 = \frac{(q - 2) (q + 1) \left(q^5 + q^4 - 6 q^3 + 3 q + 2\right)}{q^3 - q - 1}. \]
	We have
	\[ B_4 := \left| \left(q^2 - q - 2 \right)^2 - \Delta_{\Gamma_{GL(2, q)}} \right| = \left| \frac{(q - 2) (q + 1) \left(2 q^4 - 5 q^3 + 2 q + 2\right)}{q^3 - q - 1} \right|. \]
	Since $q \geq 3$ we get $(q - 2) (q + 1) (2 q^4 - 5 q^3 + 2 q + 2)= (q - 2) (q + 1)(2 + 2q + q^3(2q - 5)) > 0$. Therefore
	\[ B_4 = \frac{(q - 2) (q + 1) \left(2 q^4 - 5 q^3 + 2 q + 2\right)}{q^3 - q - 1}.\]
	We have
	\[ B_5 := \left| 2 \left(q^2 - 2 q\right) \left(q^2 - 2 q - 1\right) - \Delta_{\Gamma_{GL(2, q)}} \right| = \left| \frac{(q - 2) q \left((q - 3) (q + 1) q^3 + 5 q + 1\right)}{q^3 - q - 1} \right|. \]
	Since $q \geq 3$ we get $(q - 2) q ((q - 3) (q + 1) q^3 + 5 q + 1) > 0$ and $q^3 - q - 1 > 0$. Therefore, 
	\[ B_5 = \frac{(q - 2) q \left((q - 3) (q + 1) q^3 + 5 q + 1\right)}{q^3 - q - 1}. \]
	Also
	\[ B_6 := \left| \left(q^2 - 2 q - 1\right)^2 - \Delta_{\Gamma_{GL(2, q)}} \right| = \left| \frac{q (q (q (q (3 - 2 q) + 6) - 5) - 3) - 1}{q^3 - q - 1} \right|. \]
	Since $q \geq 3$ we get $q (q (q (q (3 - 2 q) + 6) - 5) - 3) - 1 < 0$ and $q^3 - q - 1 > 0$. Therefore
	\[ B_6 = - \frac{q (q (q (q (3 - 2 q) + 6) - 5) - 3) - 1}{q^3 - q - 1}.\]
	Hence, by \eqref{LE+cn}, we get
	\begin{align*}
		LE^+_{CN}(\Gamma_{GL(2, q)}) &= \frac{1}{2} (q^2 + 1) B_1 + \frac{1}{2} (q^2 + 1) \left(q^2 - 3 q + 1\right) B_2 + \frac{1}{2}  (q^2 - q) B_3\\
		& \quad + \frac{1}{2} (q^2 - q) (q^2 - q - 1) B_4 + (q+1) B_5 + (q+1) (q^2-2 q) B_6\\
		&= \begin{cases}
			\frac{(q - 2) q (q + 1) (q (q (q (q (2 q ((q - 4) q + 6) - 7) - 7) - 1) + 11) + 2)}{q^3 - q - 1}, &\text{ if } q\leq 5\\
			\frac{q (q + 1) \left((q (q (q (q (2 q - 11) + 20) - 16) + 7) + 8) q^3 - 12 q - 4\right)}{q^3 - q - 1}, &\text{ if } q\geq 6.
		\end{cases}
	\end{align*}

\end{proof}

\begin{prop}\label{FA-proposition}
	Let $F = GF(2^n)$ (where $n\geq 2$) and $\nu$ be the Frobenius automorphism of $F$, that is $\nu(x)=x^2$, for all $x\in F$. Then the CNL-spectrum, CNSL-spectrum, CNL-energy and CNSL-energy of the commuting graphs of the groups
	\[A(n,\nu)=\left\lbrace U(a,b)=\begin{pmatrix}
		1 & 0 & 0\\
		a & 1 & 0\\
		b & \nu(a) & 1
	\end{pmatrix}~:~a,b\in F\right\rbrace\]
	under matrix multiplication are given by\\
	\[\cnlspec(\Gamma_{A(n,\nu)})=\left\lbrace (0)^{(2^n - 1)}, (2^n(2^n - 2))^{(2^n - 1)^2}\right\rbrace,\]
	\[\cnqspec(\Gamma_{A(n,\nu)})=\left\lbrace (2(2^n - 1)(2^n - 2))^{(2^n - 1)}, ((2^n - 2)^2)^{(2^n - 1)^2}\right\rbrace\]
	and
	\[LE_{CN}(\Gamma_{A(n,\nu)}) = LE^+_{CN}(\Gamma_{A(n,\nu)}) = 2 \left(2^n - 2\right) \left(2^n - 1\right)^2.\]
\end{prop}
\begin{proof}
	By \cite[Proposition 2.4]{JR}, we have $\Gamma_{A(n,\nu)} = (2^n - 1)K_{2^n}$.
	Therefore, by Theorem \ref{CN-LE-CNSL-LE^+_Kn}, we get
	
	\centerline{	$\cnlspec(\Gamma_{A(n,\nu)})=\left\lbrace (0)^{(2^n - 1)}, (2^n(2^n - 2))^{(2^n - 1)^2}\right\rbrace$ }
	
	\noindent	and
	
	\centerline{	$\cnqspec(\Gamma_{A(n,\nu)})=\left\lbrace (2(2^n - 1)(2^n - 2))^{(2^n - 1)}, ((2^n - 2)^2)^{(2^n - 1)^2}\right\rbrace.$} 
	
	\noindent	Hence, we get the required $\cnlspec(\Gamma_{A(n,\nu)})$ and  $\cnqspec(\Gamma_{A(n,\nu)})$ on simplification.
	\par Here $|V(\Gamma_{A(n,\nu)})|= 4^n - 2^n$ and $\tr(\cnrs(\Gamma_{A(n,\nu)}))= 2^n \left(2^n-2\right) \left(2^n-1\right)^2$. So, $\Delta_{\Gamma_{A(n,\nu)}} = 4^n - 3\times 2^n + 2$.
	
	\noindent \textbf{CNL-energy calculation:}
	
We have
\[ 
L_1 := \left| 0 - \Delta_{\Gamma_{A(n,\nu)}} \right| = \left| - \left(4^n - 3\times 2^n + 2\right) \right| = 4^n - 3\times 2^n + 2, 
\]
	since $ - \left(4^n - 3\times 2^n + 2\right) \leq 0$, as $n \geq 2$.
	Also
	\[ L_2 := \left| 2^n (2^n - 2) - \Delta_{\Gamma_{A(n,\nu)}} \right| = \left| 2^n-2 \right| = 2^n - 2, \]
	since $2^n - 2 > 0$, as $n \geq 2$. Hence, by \eqref{LEcn}, we get
	\begin{align*}
		LE_{CN}(\Gamma_{A(n,p)}) &= (2^n - 1) L_1 + (2^n - 1)^2 L_2= 2 (2^n - 2) (2^n - 1)^2.
	\end{align*}
	\noindent \textbf{CNSL-energy calculation:}
	
	We have
	\[ B_1 := \left| 2 (2^n - 1) (2^n - 2) - \Delta_{\Gamma_{A(n,\nu)}} \right| = \left| 4^n - 3\times 2^n + 2 \right| = 4^n - 3\times 2^n + 2, \]
	since $ 4^n - 3\times 2^n + 2 \geq 2$, as $n \geq 2$.
	Also
	\[ B_2 := \left| (2^n - 2)^2 - \Delta_{\Gamma_{A(n,\nu)}} \right| = \left| 2 - 2^n \right| = 2^n - 2,\] 
	since $2 - 2^n < 0$, as $n\geq 2$. Hence, by \eqref{LE+cn}, we get
	\begin{align*}
		LE^+_{CN}(\Gamma_{A(n,p)}) &= (2^n - 1) B_1 + (2^n - 1)^2 B_2 = 2 \left(2^n-2\right) \left(2^n-1\right)^2.
	\end{align*}
	Hence the result follows.
\end{proof}
\begin{prop}\label{GF-proposition}
	Let $F=GF(p^n)$, where $p$ is a prime. Then the CNL-spectrum, CNSL-spectrum, CNL-energy and CNSL-energy of commuting graphs of the groups
	\[A(n,p)=\left\lbrace V(a,b,c)=\begin{pmatrix}
		1 & 0 & 0\\
		a & 1 & 0\\
		b & c & 1
	\end{pmatrix}~:~a,b,c\in F\right\rbrace\]
	under matrix multiplication $V(a,b,c)V(a',b',c')=V(a + a', b + b' + ca', c + c')$ are given by
	
	\centerline{	$\cnlspec(\Gamma_{A(n,p)})=\left\lbrace (0)^{p^n + 1}, ((p^{2n} - p^n)(p^{2n} - p^n - 2))^{(p^n + 1)(p^{2n} - p^n - 1)}\right\rbrace$,}
	
	$\cnqspec(\Gamma_{A(n,p)})$
	
	\qquad \quad	$=\left\lbrace (2(p^{2n} - p^n - 1)(p^{2n} - p^n - 2))^{(p^n + 1)}, ((p^{2n} - p^n - 2)^2)^{(p^n + 1)(p^{2n} - p^n - 1)}\right\rbrace$ 
	
	\noindent and
	$LE_{CN}(\Gamma_{A(n,p)}) = LE^+_{CN}(\Gamma_{A(n,p)}) = 2 \left(p^n + 1\right)^2 \left(-3 p^{2 n} + p^{3 n} + p^n + 2\right) $.
\end{prop}
\begin{proof}
	By \cite[Proposition 2.5]{JR},
	$\Gamma_{A(n,p)}=(p^n + 1)K_{p^{2n} - p^n}.$
	Therefore, by Theorem \ref{CN-LE-CNSL-LE^+_Kn}, we get
	
	\centerline{	$\cnlspec(\Gamma_{A(n,p)})=\left\lbrace (0)^{p^n + 1}, ((p^{2n} - p^n)(p^{2n} - p^n - 2))^{(p^n + 1)(p^{2n} - p^n - 1)}\right\rbrace$}
	
	\noindent and
	
	$\cnqspec(\Gamma_{A(n,p)})$
	
	\qquad \quad	$=\left\lbrace (2(p^{2n} - p^n - 1)(p^{2n} - p^n - 2))^{(p^n + 1)}, ((p^{2n} - p^n - 2)^2)^{(p^n + 1)(p^{2n} - p^n - 1)}\right\rbrace$. 	
	
	\par Here $|V(\Gamma_{A(n,p)})|= p^n \left(p^{2 n}-1\right)$ and  $\tr(\cnrs(\Gamma_{A(n,p)}))= p^n \left(p^n-2\right) \left(p^n-1\right) \times$ $\left(p^n+1\right)^2 \left(p^n \left(p^n-1\right)-1\right)$.  So, $\Delta_{\Gamma_{A(n,p)}} = \left(p^n-2\right) \left(p^n+1\right) \left(p^n \left(p^n-1\right)-1\right)$.
	
	\noindent \textbf{CNL-energy calculation: }
	
	We have
	\begin{align*}
		L_1 := \left| 0 - \Delta_{\Gamma_{A(n,p)}} \right| &= \left| -\left(\left(p^n-2\right) \left(p^n+1\right) \left(p^n \left(p^n-1\right)-1\right)\right) \right|\\
		&= \left(p^n-2\right) \left(p^n+1\right) \left(p^n \left(p^n-1\right)-1\right),
	\end{align*}

\noindent	since $ -\left(\left(p^n-2\right) \left(p^n+1\right) \left(p^n \left(p^n-1\right)-1\right)\right) < 0$, as $p^n-2\geq 2$, $p^n+1 > 0$ and $p^n( p^n -1 ) - 1 \geq 0.$ Also

	\[ L_2 := \left| (p^{2 n} - p^n) (p^{2 n}-p^n - 2) - \Delta_{\Gamma_{A(n,p)}} \right| = \left| (p^n - 2) (p^n + 1) \right| = (p^n - 2) (p^n + 1), \]
	since $p^n + 1 > 0$ and $p^n - 2\geq 0$, so $(p^n + 1)(p^n - 2)\geq 0$. Hence, by \eqref{LEcn}, we get
	\begin{align*}
		LE_{CN}(\Gamma_{A(n,p)}) &= (p^n + 1) L_1 + (p^{2 n} - p^n - 1) (p^n + 1) L_2\\
		&= 2 (p^n + 1)^2 (- 3 p^{2 n} + p^{3 n} + p^n + 2).
	\end{align*}

	\noindent \textbf{CNSL-energy calculation:}
	
	We have
	\begin{align*}
		B_1 := &\left| 2 (p^{2 n} - p^n - 1) (p^{2 n} - p^n - 2) - \Delta_{\Gamma_{A(n,p)}} \right|\\ 
		= &\left| \left(p^n-2\right) \left(p^n+1\right) \left(p^n \left(p^n-1\right)-1\right) \right|= \left(p^n-2\right) \left(p^n+1\right) \left(p^n \left(p^n-1\right)-1\right),
	\end{align*}
	since $\left(p^n-2\right) \left(p^n+1\right) \left(p^n \left(p^n-1\right)-1\right) > 0$, as $(p^n - 2) \geq 0$, $p^n + 1\geq 3$, $p^n - 1 \geq 1$ and $p^n(p^n - 1) - 1 > 0$. 
	Also
	\[B_2 := \left| (p^{2 n} - p^n - 2)^2 - \Delta_{\Gamma_{A(n,p)}} \right| = \left| - p^{2 n} + p^n + 2 \right| = p^{2 n} - p^n - 2,\]
	since $ - p^{2 n} + p^n + 2 = p^n(1 - p^n) + 2 \leq 0$. Hence, by \eqref{LE+cn}, we get
	\begin{align*}
		LE^+_{CN}(\Gamma_{A(n,p)}) &= (p^n + 1) B_1 + \left(p^{2 n} - p^n - 1\right) \left(p^n + 1\right) B_2\\
		&= 2 (p^n + 1)^2 ( - 3 p^{2 n} + p^{3 n} + p^n + 2).
	\end{align*}
	Hence the result follows.
\end{proof}

\subsection{Groups whose central quotient is isomorphic to $Sz(2)$, $D_{2m}$ or $\mathbb{Z}_p \times \mathbb{Z}_p$ }
Let us begin with the groups $G$ such that $G/Z(G)$ is  isomorphic to the Suzuki group $Sz(2)=\langle a,b~:~a^5 = b^4 = 1,~b^{-1}ab = a^2 \rangle$.
\begin{thm}\label{Suzuki-theorem}
	Let $G$ be a finite non-abelian group such that $\frac{G}{Z(G)}\cong Sz(2)$ and $|Z(G)| = z$. Then the CNL-spectrum, CNSL-spectrum, CNL-energy and CNSL-energy of $\Gamma_G$ are given by
	$$\cnlspec(\Gamma_G) = \left\lbrace (0)^6, (4z(4z - 2))^{(4z - 1)}, (3z(3z - 2))^{5(3z - 1)}\right\rbrace,$$	
	$\cnqspec(\Gamma_G)$
	
	 $= \left\lbrace(2(4z - 1)(4z - 2))^1, ((4z - 2)^2)^{(4z - 1)}, (2(3z - 1)(3z - 2))^5, ((3z - 2)^2)^{5(3z - 1)} \right\rbrace, $
	\[LE_{CN}(\Gamma_G) = \begin{cases}
								\frac{648}{19}, &\text{ for } z=1\\
								\frac{2}{19} (4 z - 1) (z (105 z + 31) - 38), &\text{ for } z\geq 2
						 \end{cases}
	\]
and
	\[LE^+_{CN}(\Gamma_G) = \frac{10}{19} (3 z - 1) (z (28 z + 45) - 38).\]
\end{thm}
\begin{proof}
	By \cite[Theorem 2.2]{JR}, we have
	$\Gamma_G=K_{4z}\sqcup 5 K_{3z}.$
	Therefore, by Theorem \ref{CN-LE-CNSL-LE^+_Kn}, we get
	
	$$\cnlspec(\Gamma_G)=\left\lbrace (0)^1, (4z(4z - 2))^{(4z - 1)}, 0^5, (3z(3z - 2))^{5(3z - 1)}\right\rbrace$$ 
	and 
	$\cnqspec(\Gamma_G)=\left\lbrace(2(4z - 1)(4z - 2))^1, ((4z - 2)^2)^{(4z - 1)}, (2(3z - 1)(3z - 2))^5,\right. \\ \left. ((3z - 2)^2)^{5(3z - 1)} \right\rbrace. $ Hence we get the required $\cnlspec(\Gamma_G)$ and $\cnqspec(\Gamma_G)$ on simplification.
	\par Here $|V(\Gamma_G)|= 19z$ and $\tr(\cnrs(\Gamma_G))= 199 z^3 - 183 z^2 + 38 z$. So, $ \Delta_{\Gamma_G} = \frac{199 z^2 - 183 z + 38}{19}$. 
	
	\noindent \textbf{CNL-energy calculation:}
	
	We have
	\[ L_1 := \left| 0 - \Delta_{\Gamma_G} \right| = \left| - \frac{199 z^2 - 183 z + 38}{19} \right| = \frac{199 z^2 - 183 z + 38}{19}, \]
	since $199 z^2 - 183 z + 38  = z(199 z - 183)  + 38 > 0$, for all $z \geq 1$.\\
	We have
	\[ L_2 := \left| 4 z (4 z-2) - \Delta_{\Gamma_G} \right| = \left| \frac{105 z^2 + 31 z - 38}{19} \right| = \frac{105 z^2 + 31 z - 38}{19},\] 
	since, $105 z^2 + 31 z - 38 > 0$, for all $z\geq 1$. \\
	Also
	\[ L_3 := \left|3 z (3 z-2)- \Delta_{\Gamma_G}\right| = \left|  \frac{z(69 - 28 z) - 38}{19} \right|. \]
	For $z = 1$, $z(69 - 28 z) - 38= 3 > 0$ and for $z\geq 2$, $z(69 - 28 z) - 38 < 0$. Therefore
	\[
	L_3 = \begin{cases}
		\frac{3}{19}, &\text{ for } z=1\\
		- \frac{z(69 - 28 z) - 38}{19}, &\text{ for } z\geq 2.
	\end{cases}
	\]
Hence, by \eqref{LEcn}, we get
\begin{align*}
	LE_{CN}(\Gamma_G) &= 6 L_1 + (4z - 1) L_2 + 5(3z - 1)L_3\\
					  &= \begin{cases}
						\frac{648}{19}, &\text{ for } z = 1\\
						\frac{2}{19} (4 z - 1) (z (105 z + 31) - 38), &\text{ for } z\geq 2.
					\end{cases}
\end{align*}

\noindent\textbf{CNSL-energy calculation:}

We have
\[ B_1 := \left| 2 (4 z - 1) (4 z - 2) - \Delta_{\Gamma_G} \right| = \left| \frac{409 z^2 - 273 z + 38}{ 19 } \right| = \frac{409 z^2 - 273 z + 38}{ 19 }, \]
since $409 z^2 - 273 z + 38 = ( 409 z - 273 ) z + 38 > 0$, as $z\geq 1$.\\
We have
\[ B_2 := \left| (4 z - 2)^2 - \Delta_{\Gamma_G} \right| = \left| \frac{105 z^2 - 121 z + 38}{ 19 } \right| = \frac{105 z^2 - 121 z + 38}{ 19 }, \]
since $105 z^2 - 121 z + 38 = ( 105 z - 121 ) z + 38 > 0$, for all $k\geq 1$.

We have

\[ B_3 := \left| 2 (3 z-1) (3 z-2) - \Delta_{\Gamma_G} \right| = \left| \frac{143 z^2 - 159 z + 38}{ 19 } \right| = \frac{143 z^2 - 159 z + 38}{ 19 }, \]
since $143 z^2 - 159 z + 38 = ( 143 z - 159 ) z + 38 > 0$, for all $z\geq 1$.

Also

\[ B_4 := \left| (3 z-2)^2 - \Delta_{\Gamma_G} \right| = \left| \frac{-28 z^2 - 45 z + 38}{ 19 } \right| = \frac{28 z^2 + 45 z - 38}{ 19 }, \]
since $ - 28 z^2 - 45 z + 38 = - z (28 z + 45) + 38 < 0$, for all $z\geq 1$. 
Therefore, by \eqref{LE+cn}, we get

\begin{align*}
	LE^+_{CN}(\Gamma_G) &= 1 \times B_1 + (4z - 1) B_2 + 5 \times B_3 + 5(3z - 1) B_4\\
						&= \frac{10}{19} (3 z-1) (z (28 z+45)-38).
\end{align*}
	Hence the result follows.
\end{proof}

\begin{thm}\label{Z_p-theorem}
	Let $G$ be a finite non-abelian group such that $\frac{G}{Z(G)}\cong \mathbb{Z}_p\times \mathbb{Z}_p$, where $p$ is a prime. Then the CNL-spectrum, CNSL-spectrum, CNL-energy and CNSL-energy of the commuting graph of $G$ are given by
	$$\cnlspec(\Gamma_G) = \left\lbrace (0)^{(p + 1)}, (((p - 1)z)((p - 1)z - 2))^{(p + 1)((p - 1)z - 1)}\right\rbrace,$$
	$\cnqspec(\Gamma_G)\!=\!\left\lbrace (2((p -\!1)z -\!1)((p - 1)z - 2))^{(p + 1)}, (((p - 1)z - 2)^2)^{(p + 1)((p - 1)z - 1)}\right\rbrace$ and\\
	$LE_{CN}(\Gamma_G) = LE^+_{CN}(\Gamma_G)=\begin{cases}
		0,~~~~~&\text{ if } z = 1,~p = 2\\
		2 (p + 1) ((p - 1) z - 2) ((p - 1) z - 1),~~~~~&\text{ otherwise; }
	\end{cases}$ \\
where $z = |Z(G)|.$	
\end{thm}
\begin{proof}
From \cite[Theorem 2.1]{JR-16}, we have
$\Gamma_G=(p + 1)K_{(p - 1)z}.$	
Therefore, by Theorem \ref{CN-LE-CNSL-LE^+_Kn}, we get
$$\cnlspec(\Gamma_G)=\left\lbrace (0)^{(p + 1)}, (((p - 1)z)((p - 1)z - 2))^{(p + 1)((p - 1)z - 1)}\right\rbrace$$ and 
\begin{align*}
&\cnqspec(\Gamma_G)\\
&\qquad =\left\lbrace (2((p - 1)z - 1)((p - 1)z - 2))^{(p + 1)}, (((p - 1)z - 2)^2)^{(p + 1)((p - 1)z - 1)}\right\rbrace. 
\end{align*}
Hence we get the required $\cnlspec(\Gamma_G)$ and $\cnqspec(\Gamma_G)$ on simplification.
\par Here $|V(\Gamma_G)|= \left(p^2 - 1\right) z$ and $\tr(\cnrs(\Gamma_G))= (p-1) (p+1) z ((p-1) z-2) ((p-1) z-1)$. So, $\Delta_{\Gamma_G} =  ((p - 1) z - 2) ((p - 1) z - 1)$. 

\vspace{.3cm}

\noindent \textbf{CNL-energy calculation:}

We have
\begin{align*}
	L_1 := \left| 0 - \Delta_{\Gamma_G} \right| &= \left| - (((p - 1) z - 2) ((p - 1) z - 1)) \right|\\
	&= ((p-1) z-2) ((p-1) z-1),
\end{align*}
since $((p - 1) z - 2) ((p - 1) z - 1) \geq 0$, as $p \geq 2$ and $z \geq 1$.\\
Also
\begin{align*}
	L_2 := \left| ((p - 1) z) ((p - 1) z - 2) - \Delta_{\Gamma_G} \right| &= \left| (p-1) z-2 \right|\\ &= \begin{cases}
		2 - (p - 1) z= 1, &\text{ if }z = 1,~p = 2\\
		(p - 1) z - 2, &\text{ otherwise }	
	\end{cases}	
\end{align*}
as for $z \geq 2$, $(p - 1) z - 2 > 0$ and for $z = 1$, $p = 2$, $(p - 1) z - 2 = - 1 < 0$.
Hence, by \eqref{LEcn}, we get
\begin{align*}
	LE_{CN}(\Gamma_G) &= (p + 1) L_1 + (p + 1)((p - 1)z - 1) L_2\\
					  &= \begin{cases}
						0,~~~~~&\text{ if } z = 1, p = 2\\
						2 (p + 1) ((p - 1) z - 2) ((p - 1) z - 1), &\text{ otherwise. }
						\end{cases}
\end{align*}
\noindent\textbf{CNSL-energy calculation:}

We have
\begin{align*}
	B_1 :=\left| 2 ((p - 1) z - 1) ((p - 1) z - 2) - \Delta_{\Gamma_G} \right| &= \left| ((p-1) z-2) ((p-1) z-1) \right|\\
		&= ((p-1) z-2) ((p-1) z-1),
\end{align*}
since $((p-1) z-2)((p-1) z-1) \geq 0$, as $p \geq 2$ and $z \geq 1$.\\
Also
\begin{align*}
	B_2 := \left| ((p - 1) z - 2)^2 - \Delta_{\Gamma_G} \right| &= \left| - p z + z + 2 \right|\\
	&= \begin{cases}
		2 + z - pz = 1, &\text{ if }z = 1,~p = 2\\
		pz - z - 2, &\text{ otherwise } 
	\end{cases}
\end{align*}
as for $z \geq 2$, $(p - 1) z - 2 > 0$ and for $z = 1$, $p = 2$, $(p - 1) z - 2 = - 1 < 0$. Hence, by \eqref{LE+cn}, we get
\begin{align*}
	LE^+_{CN}(\Gamma_G) &= (p + 1) B_1 + (p + 1)((p - 1)z - 1) B_2\\
						&= \begin{cases}
							0, &\text{ if }z = 1,~p = 2\\
							2(p + 1) ((p - 1)z - 1) ((p - 1)z - 2), &\text{ otherwise. }
						\end{cases}
\end{align*}
\end{proof}

\begin{cor}\label{Z_p-cor1}
	Let $G$ be a non-abelian group of order $p^3$, where $p$ is a prime. Then the CNL-spectrum, CNSL-spectrum, CNL-energy and CNSL-energy of the commuting graph of $G$ are given by

	\centerline{$\cnlspec(\Gamma_G)=\left\lbrace (0)^{(p + 1)}, (p(p - 1)(p(p - 1) - 2))^{(p + 1)((p - 1)p - 1)}\right\rbrace$,}
\noindent $\cnqspec(\Gamma_G)\!=\!\left\lbrace(2((p -\!1)p -\!1)((p - 1)p - 2))^{(p + 1)}, (((p - 1)p - 2)^2)^{(p + 1)((p - 1)p - 1)}\right\rbrace$ and
	$LE_{CN}(\Gamma_G)=LE^+_{CN}(\Gamma_G)= 2 (p + 1) ((p - 1) p - 2) ((p - 1) p - 1)$.
	

\end{cor}
\begin{proof}
	We know that for a non-abelian group $G$ of order $p^3$, $|Z(G)|=p$ and $\frac{G}{Z(G)}\cong \mathbb{Z}_p\times \mathbb{Z}_p$. Hence the result follows from the Theorem \ref{Z_p-theorem}.  
\end{proof}

\begin{cor}\label{Z_p-cor2}
	Let $G$ be a finite $4$-centralizer group. Then the CNL-spectrum, CNSL-spectrum, CNL-energy and CNSL-energy of the commuting graph of $G$ are given by
	
	\centerline{$\cnlspec(\Gamma_G) = \left\lbrace (0)^3, ((z(z - 2))^{3((p - 1)z - 1)}\right\rbrace$,}
	
	\centerline{$\cnqspec(\Gamma_G) = \left\lbrace (2(z - 1)(z - 2))^3, ((z - 2)^2)^{3(z - 1)}\right\rbrace$ and} 
	
\centerline{	$LE_{CN}(\Gamma_G)= LE^+_{CN}(\Gamma_G)= \begin{cases}
		0,~~~~~~~~ \text{ if } p=2, z=1\\
		6 (z-2) (z-1),~~~~~~~~ \text{ otherwise; }
	\end{cases}$}

\noindent where $z = |Z(G)|.$ 
\end{cor}
\begin{proof}
	Here $\frac{G}{Z(G)}=\mathbb{Z}_2\times \mathbb{Z}_2$. Hence the results follows from the Theorem \ref{Z_p-theorem}.
\end{proof} 
\begin{cor}\label{Z_p-cor3}
	Let $G$ be a finite $5$-centralizer group. Then the CNL-spectrum, CNSL-spectrum, CNL-energy and CNSL-energy of the commuting graph of $G$ are given by
	
\centerline{	$\cnlspec(\Gamma_G)=\left\lbrace (0)^4, (2z(2z - 2))^{4(2z - 1)}\right\rbrace$, }

\centerline{	$\cnqspec(\Gamma_G)=\left\lbrace (2(2z - 1)(2z - 2))^4, ((2z - 2)^2)^{4(z - 1)}\right\rbrace$}

 and 
	$LE_{CN}(\Gamma_G)= LE^+_{CN}(\Gamma_G)= 8 (2 z-2) (2 z-1)$,
	where $z = |Z(G)|.$ 

\end{cor}
\begin{proof}
	We have $\frac{G}{Z(G)}\cong \mathbb{Z}_3\times \mathbb{Z}_3$. Hence the result follows from Theorem \ref{Z_p-theorem}.
\end{proof}

\begin{cor}
Let  $G$ be a finite $(p+2)$-centralizer $p$-group. Then the CNL-spectrum, CNSL-spectrum, CNL-energy and CNSL-energy of the commuting graph of $G$ are given by
$$\cnlspec(\Gamma_G) = \left\lbrace (0)^{(p + 1)}, (((p - 1))((p - 1)|Z(G)| - 2))^{(p + 1)((p - 1)|Z(G)| - 1)}\right\rbrace,$$
$\cnqspec(\Gamma_G) = \left\lbrace (2((p -\!1)|Z(G)| -\!1)((p - 1)|Z(G)| - 2))^{(p + 1)}\right.,$

\qquad\qquad\qquad\qquad\qquad\qquad\qquad $\left. (((p - 1)|Z(G)| - 2)^2)^{(p + 1)((p - 1)|Z(G)| - 1)}\right\rbrace$ and\\
$LE_{CN}(\Gamma_G) = LE^+_{CN}(\Gamma_G)=
	2 (p + 1) ((p - 1) |Z(G)| - 2) ((p - 1)|Z(G)| - 1).$ 
\end{cor}
\begin{proof}
	We have $\frac{G}{Z(G)}\cong \mathbb{Z}_p\times \mathbb{Z}_p$. Hence the result follows from Theorem \ref{Z_p-theorem}.
\end{proof}

\begin{thm}\label{D_2m theorem}
	If $G$ is a finite group such that $\frac{G}{Z(G)}\cong D_{2m}$, $m\geq 2$, then the CNL-spectrum, CNSL-spectrum, CNL-energy and CNSL-energy of the commuting graph of $G$ are given by
	
\centerline{	$\cnlspec(\Gamma_G)=\left\lbrace (0)^{m+1}, ((m - 1)z((m - 1)z - 2))^{((m - 1)z - 1)}, (z(z-2))^{m(z-1)}\right\rbrace,$}

	$\cnqspec(\Gamma_G) = \left\lbrace (2((m - 1)z - 1)((m - 1)z - 2))^1, (((m - 1)z - 2)^2)^{((m - 1)z - 1)},\right.$
	
	$\qquad\qquad\qquad\qquad\qquad\qquad\qquad\qquad\qquad\quad\,\left. (2(z-1)(z-2))^{m}, ((z-2)^2)^{m(z-1)}\right\rbrace,$
	\[LE_{CN}(\Gamma_G) = \begin{cases}
							0, &\text{ for } m=2, z=1\\
							\frac{-6 \left(m^3+1\right) z+4 \left(2 m^2+m-1\right)+2 \left(m \left(m (m-1)^2+3\right)-1\right) z^2}{2 m-1}, &\text{ for } m=2, z \geq 2\\
							\frac{2 ((m-1) z-1) (m (z ((m-2) m z-m+3)-4)+z+2)}{2 m-1}, &\text{ otherwise}
						\end{cases}\]
	\[\text{and } LE^+_{CN}(\Gamma_G) = \begin{cases}
							0, &\text{ for } m=2, z=1\\
							\frac{1}{3} \left(18 z^2-54 z+36\right), &\text{ for } m=2, z \geq 2\\
							\frac{2 m \left(\left((m-3) m^2+1\right) z^2-3 ((m-5) m+3) z-4 m+2\right)}{2 m-1},&\text{ for } m=3, z=1 \\
							&\text{ and } m=4, z=1\\
							\frac{2 (m-2) (m-1) m z^2 (m z-3)}{2 m-1}, &\text{ otherwise;}
\end{cases}\]
where $z=|Z(G)|$.

\end{thm}

\begin{proof}
	From \cite[Theorem 2.5]{JR-16}, we have
	$\Gamma(G)=K_{(m - 1)z}\sqcup m K_{z}.$
	Therefore, by Theorem \ref{CN-LE-CNSL-LE^+_Kn}, we get\\
	$\cnlspec(\Gamma_G) = \left\lbrace (0)^{m + 1}, ((m - 1)z((m - 1)z - 2))^{((m - 1)z - 1)}, (z(z - 2))^{m(z - 1)}\right\rbrace$ and \\
	$\cnqspec(\Gamma_G)= \left\lbrace (2((m - 1)z - 1)((m - 1)z - 2))^1, (((m - 1)z - 2)^2)^{((m - 1)z - 1)},\right.$
	
	$\qquad\qquad\qquad\qquad\qquad\qquad\qquad\qquad\qquad\quad\,\left. (2(z - 1)(z - 2))^{m}, ((z - 2)^2)^{m(z - 1)}\right\rbrace$.
	
\noindent	Hence we get the required $\cnlspec(\Gamma_G)$ and $\cnqspec(\Gamma_G)$ on simplification.
	\par Here $|V(\Gamma_G)|= (2 m - 1) z$ and $\tr(\cnrs(\Gamma_G))= m (z - 2) (z - 1) z + (m - 1) z ((m - 1) z - 2) ((m - 1) z - 1)$. So, $\Delta_{\Gamma_G} = \frac{(m ((m - 3) m + 4) - 1) z^2 - 3 ((m - 1) m + 1) z + 4 m - 2}{2 m - 1}$.

	\noindent \textbf{CNL-energy calculation:}
	
	We have
	\[ L_1 \! :=\! \left| 0\! -\! \Delta_{\Gamma_G} \right|\! =\! \left| \frac{(2 - 4m) \!+\! 3z(1 - m)\! +\! z^2(1 - 4m) + m^2 z^2(3 - \frac{m}{2}) + m^2 z(3 - \frac{m z}{2})}{2m - 1} \right|.\]
Let $\alpha_1(m,z)= (2 - 4m) + 3z(1 - m) + z^2(1 - 4m) + m^2 z^2(3 - \frac{m}{2}) + m^2 z(3 - \frac{m z}{2})$. For $m \geq 7$, $2 - 4m < 0$, $1 - 3m < 0$, $1 - 4m < 0$, $3 - \frac{m}{2} < 0$ and $3 - \frac{m z}{2} < 0$ as $z \geq 1$, so $\alpha_1(m,z) < 0$. Also, $\alpha_1(2,z)= - 3 z^2 + 9 z - 6 \leq 0$, $\alpha_1(3,z)= - 11 z^2 + 21 z - 10 \leq 0$, $\alpha_1(4,z)= - 31 z^2 + 39 z - 14 \leq 0$, $\alpha_1(5,z)= - 69 z^2 + 63 z - 18 \leq 0$ and $\alpha_1(6,z)=- 131  z^2 + 93 z - 22 \leq 0$ as $z\geq 1$. Therefore,
\[L_1 = - \frac{(2 - 4m) + 3z(1 - 3m) + z^2(1 - 4m) + m^2 z^2(3 - \frac{m}{2}) + m^2 z(3 - \frac{m z}{2})}{2m - 1}.\] 
We have
\[ L_2 := \left| ((m - 1) z) ((m - 1) z - 2) - \Delta_{\Gamma_G} \right| = \left| \frac{m (z ((m - 2) m z - m + 3) - 4) + z + 2}{2 m - 1} \right|. \]
Let $\alpha_2(m,z) = m (z ((m - 2) m z - m + 3) - 4) + z + 2$. For $m \geq 4$, $3 - m\geq 1$, $m - 2\geq 2$ so $m(m - 2)\geq 8$. Therefore, $(m - 2) m z - m + 3 \geq 7$ which gives $m (z ((m - 2) m z - m + 3) - 4) > 0$ and so $\alpha_2(m,z) > 0$. Also, $\alpha_2(2,1)= - 3$ and $\alpha_2(2,z)= 3z - 6 \geq 0$ for all $z \geq 2$. Further, $\alpha_2(3,z) = 9z^2 + z - 10 \geq 0$ for all $z \geq 1$. Hence,
\[L_2 = \begin{cases}
	- \frac{m (z ((m - 2) m z - m + 3) - 4) + z + 2}{2 m - 1}, & \text{ for } m=2,~ z=1\\
	\frac{m (z ((m - 2) m z - m + 3) - 4) + z + 2}{2 m - 1}, & \text{ otherwise. }
\end{cases}\]
Also
\[L_3 := \left| z (z - 2) - \Delta_{\Gamma_G} \right| = \left| \frac{m (- z ((m - 2) (m - 1) z - 3 m + 7) - 4) + 5 z + 2}{2 m - 1} \right|.\]
Let $\alpha_3(m,z)= m (- z ((m - 2) (m - 1) z - 3 m + 7) - 4) + 5 z + 2 = 2 - 4m + 5z - 7mz + (mz(3m - (m - 2)(m - 1)z))$. For $m \geq 6$, $3m - (m - 2)(m - 1) = - m^2 + 6m - 2 < 0$ which gives $3m <  (m - 2)(m - 1)$ and so $3m - (m - 2)(m - 1)z < 0$.  Thus, $mz (3m - (m - 2)(m - 1)z) <0$. Therefore, $\alpha_3(m,z) < 0$, for all $m \geq 6$, as $2 - 4m < 0$ and $5z - 7mz < 0$. We have $\alpha_3(2,z) = 3z - 6$ which is $< 0$ if $z=1$ and $\geq 0$ if $z\geq 2$.
Also, $\alpha_3(3,z) = - 6z^2 + 11z -10 \leq 0$, $\alpha_3(4,z) = - 24z^2 + 25z - 14 \leq 0$ and $\alpha_3(5,z) = - 60z^2 + 45z - 18 \leq 0$ for all $z \geq 1$. Hence
\[L_3 = \begin{cases}
		\frac{m (-z ((m-2) (m-1) z-3 m+7)-4)+5 z+2}{2 m-1}, &\text{ for } m = 2, z\geq 2\\
		- \frac{m (-z ((m-2) (m-1) z-3 m+7)-4)+5 z+2}{2 m-1}, &\text{ otherwise}.
\end{cases} \]
Therefore, by \eqref{LEcn}, we get
\begin{align*}
	LE_{CN}(\Gamma_G) &= (m+1) L_1 + ((m - 1)z - 1) L_2 + m(z - 1) L_3\\
					  &= \begin{cases}
					  	0, &\text{ for } m=2, z=1\\
					  	\frac{-6 \left(m^3+1\right) z+4 \left(2 m^2+m-1\right)+2 \left(m \left(m (m-1)^2+3\right)-1\right) z^2}{2 m-1}, &\text{ for } m=2, z \geq 2\\
					  	\frac{2 ((m-1) z-1) (m (z ((m-2) m z-m+3)-4)+z+2)}{2 m-1}, &\text{ otherwise}.
					  \end{cases}
\end{align*}
	\noindent \textbf{CNSL-energy calculation:}
	
	We have
	\begin{align*}
		B_1 &:= \left| 2 ((m - 1) z - 1) ((m - 1) z - 2) - \Delta_{\Gamma_G} \right|\\
		&= \left| \frac{((m - 1) m (3 m - 4) - 1) z^2 - 3 (m (3 m - 5) + 1) z + 4 m - 2}{2 m - 1} \right|.
	\end{align*}
Let $\beta_1(m,z) = ((m - 1) m (3 m - 4) - 1) z^2 - 3 (m (3 m - 5) + 1) z + 4 m - 2$. Note that $((m - 1)m(3m - 4) -1) - 3(m(3m - 5) + 1) = (3m - 4)(m(m - 4) + 1) > 0$ for all $m \geq 4$. Therefore, for all $m \geq 4$ and $z \geq 1$, $((m - 1) m (3 m - 4) - 1) z^2 - 3 (m (3 m - 5) + 1) z > 0$ and so $\beta_1(m,z) > 0$ since $4m - 2 > 0$. Also, $\beta_1(2,z) = 3z^2 - 9z + 6\geq 0$ and $\beta_1(3,z) = 29z^2 - 39z + 10 \geq 0$ for all $z \geq 1$. Hence
\[B_1 = \frac{((m - 1) m (3 m - 4) - 1) z^2 - 3 (m (3 m - 5) + 1) z + 4 m - 2}{2 m - 1}. \]
We have
\[B_2 := \left| ((m-1) z-2)^2 - \Delta_{\Gamma_G} \right| = \left| \frac{m (z ((m - 2) m z - 5 m + 9) + 4) - z - 2}{2 m - 1} \right|.\]
Let  $\beta_2(m,z) = m (z ((m - 2) m z - 5 m + 9) + 4) - z - 2 = 4m - 2 - z + 9mz + mz(mz(m - 2) - 5m)$. Clearly, for $m \geq 7$ and  $z \geq 1$, $\beta_2(m,z) > 0$ as $mz(m - 2) - 5m \geq 0$. Also, $\beta_2(2,z) = 6 -3z \leq 0$
for all $z\geq 2$ and $\beta_2(2,1) = 3$. Further, $\beta_2(3,z) = 9z^2 - 19z + 10 \geq 0$, $\beta_2(4,z) = 32z^2 - 45z + 14 \geq 0$, $\beta_2$, $(5,z) = 75z^2 - 81z + 18$ and $\beta_2(6,z) = 144z^2 - 127z + 22 \geq 0$ for all $z \geq 1$. Therefore
\[B_2 = \begin{cases}
		- \frac{m (z ((m - 2) m z - 5 m + 9) + 4) - z - 2}{2 m - 1}, &\text{ for } m=2,  z\geq 2\\
		\frac{m (z ((m - 2) m z - 5 m + 9) + 4) - z - 2}{2 m - 1}, &\text{ otherwise. }
\end{cases}\]
We have
\begin{align*} 
	B_3 &:= \left| 2 (z - 1) (z - 2) - \Delta_{\Gamma_G} \right| \\
	&= \left| \frac{ - \left(\left((m - 3) m^2 + 1\right) z^2\right) + 3 ((m - 5) m + 3) z + 4 m - 2}{2 m - 1} \right|.
\end{align*}
Let $\beta_3(m,z) = - \left(\left((m - 3) m^2 + 1\right) z^2\right) + 3 ((m - 5) m + 3) z + 4 m - 2 = - 2 + 2m(2 - 3z) + 9z(1 - m) + m^2 z(3 - \frac{m z}{2}) - z^2 + m^2 z^2 (3 - \frac{m}{2})$. For $m\geq 6$, $3 - \frac{m}{2} \leq 0$, $1 - m \leq 0$ and $3 - \frac{m z}{2} \leq 0$ as  $z \geq 1$. Therefore, $\beta_3(m,z) \leq 0$ for all $m\geq 6$ and  $z \geq 1$. Also, $\beta_3(2,z)= 3z^2 - 9z + 6 \geq 0$, $\beta_3(3,z) = - z^2 - 9z + 10 \leq 0$, $\beta_3(4,z)= - 17z^2 - 3z + 14 \leq 0$ and $\beta_3(5,z) = - 51z^2 + 9z + 18 \leq 0$ for all $z \geq 1$. So
\[B_3 = \begin{cases}
	 \frac{-\left(\left((m-3) m^2+1\right) z^2\right)+3 ((m-5) m+3) z+4 m-2}{2 m-1}, &\text{ for } m=2\\
	-  \frac{-\left(\left((m-3) m^2+1\right) z^2\right)+3 ((m-5) m+3) z+4 m-2}{2 m-1}, &\text{ for } m\geq 3.
\end{cases}\]
We have
\[ B_4 := \left| (z-2)^2- \Delta_{\Gamma_G} \right| = \left| \frac{m (4-z ((m-2) (m-1) z-3 m+11))+7 z-2}{2 m-1} \right|.\]
Let $\beta_4(m,z) = m (4 - z ((m - 2) (m - 1) z - 3 m + 11)) + 7 z - 2 = - 2 + 7z - 11zm + m(4 - z( - 3m + (m - 2)(m - 1)z))$. For $m \geq 7$, $(m - 2)(m - 1) - 3m = m^2 - 6m + 2 \geq 9$ which gives $(m - 2)(m - 1)z - 3m \geq 9$  (as $z\geq 1$).  Therefore, $4 - z ((m - 2)(m - 1)z - 3m) \leq - 5$ and so $m(4 - z ((m - 2)(m - 1)z - 3m))\leq -5 < 0$. Hence, $\beta_4(m,z) < 0$ for all $m \geq 7$ and $z \geq 1$, as $7z - 11zm < 0$. Also
 \\ $\beta_4(2,z) = 6 - 3z  \begin{cases}
	3 > 0, &\text{ for } z = 1\\
	\leq 0, &\text{ for } z \geq 2,
\end{cases}$ $\beta_4(3,z) = - 6z^2 + z + 10 \begin{cases}
5 > 0, &\text{ for } z = 1\\
\leq 0, &\text{ for } z \geq 2,
\end{cases}$\\ $\beta_4(4,z) = - 24z^2 + 11z + 14 \begin{cases}
1 > 0, &\text{ for } z = 1\\
\leq 0, &\text{ for } z \geq 2,
\end{cases}$
 $\beta_4(5,z) = - 60z^2 + 27z + 18 < 0$ and $\beta_4(6,z) = - 120z^2 + 49z + 22 < 0$ for all $z \geq 1$. Therefore,
\[B_4 = \begin{cases}
	\frac{m (4 - z ((m - 2) (m - 1) z - 3 m + 11)) + 7 z - 2}{2 m - 1}, &\text{for } m = 2, z = 1; m = 3, z = 1;\\
	&m = 4, z = 1\\
	- \frac{m (4 - z ((m - 2) (m - 1) z - 3 m + 11)) + 7 z - 2}{2 m - 1}, &\text{otherwise. }
\end{cases}
\]
Hence, by \eqref{LE+cn}, we get
\begin{align*}
	LE^+_{CN}(\Gamma_G) &= 1\times B_1 + ((m - 1)z - 1)\ B_2 + m\ B_3 + m(z - 1)\ B_4\\
						&= \begin{cases}
							0, &\text{ for } m=2, z=1\\
							\frac{1}{3} \left(18 z^2-54 z+36\right), &\text{ for } m=2, z \geq 2\\
							\frac{2 m \left(\left((m-3) m^2+1\right) z^2-3 ((m-5) m+3) z-4 m+2\right)}{2 m-1},&\text{ for } m=3, z=1 \\ 
							&\text{ and } m=4, z=1\\
							\frac{2 (m-2) (m-1) m z^2 (m z-3)}{2 m-1}, &\text{ otherwise}.
						\end{cases}
\end{align*}
\end{proof}

\begin{cor}\label{D_2m-cor1}
	The CNL-spectrum, CNSL-spectrum, CNL-energy and CNSL-energy of the commuting graphs of the metacyclic groups $\mathcal{M}_{2mn}=\langle x, y: x^m=y^{2n}=1,~yxy^{-1}=x^{-1} \rangle$, where $m \geq 2$, are given by\\
	$\cnlspec(\Gamma_{\mathcal{M}_{2mn}})=$
	
	$\begin{cases}
		\left\lbrace (0)^{m + 1}, ((m - 1)n((m - 1)n - 2))^{((m - 1)n - 1)},  (n(n - 2))^{m(n - 1)} \right\rbrace,  &\text{ if }2\nmid m\\
		 
		\left\lbrace (0)^{\frac{m}{2} + 1}, ((\frac{m}{2} - 1)2n((\frac{m}{2} - 1)2n - 2))^{((\frac{m}{2} - 1)2n - 1)},  (2n(2n - 2))^{\frac{m}{2}(2n - 1)}\right\rbrace,  &\text{ if }2\mid m,
	\end{cases}$

	$\cnqspec(\Gamma_{\mathcal{M}_{2mn}})=$
	
\qquad	$\begin{cases}
		\left\lbrace (2((m - 1)n - 1)((m - 1)n - 2))^1,  (((m - 1)n - 2)^2)^{((m - 1)n - 1)}, \right. \\ 
		\qquad \qquad \qquad \qquad \qquad \left. (2(n - 1)(n - 2))^{m}, ((n - 2)^2)^{m(n - 1)}\right\rbrace,  \quad \text{ if }2\nmid m\\
		
		\left\lbrace (2((\frac{m}{2} - 1)2n - 1)((\frac{m}{2} - 1)2n - 2))^1, (((\frac{m}{2} - 1)2n - 2)^2)^{((\frac{m}{2} - 1)2n - 1)}, \right. \\ 
		\qquad \qquad \qquad \qquad \left. (2(2n - 1)(2n - 2))^{\frac{m}{2}}, ((2n - 2)^2)^{\frac{m}{2}(2n -1 )}\right\rbrace,  \quad\text{ if }2\mid m,
	\end{cases}$
\[
LE_{CN}(\Gamma_{\mathcal{M}_{2mn}})= \begin{cases}
	\frac{2 ((m-1) n-1) (m (n ((m-2) m n-m+3)-4)+n+2)}{2 m-1}, &\text{ if } 2 \nmid m\\
	\frac{1}{3} \left(72 n^2-108 n+36\right), &\text{ if } m=2\\
	\frac{((m-2) n-1) (m (n ((m-4) m n-m+6)-4)+4 (n+1))}{m - 1}, &\text{ if } m\geq 4 \text{ is even}
\end{cases}
\]
and
\[
LE^+_{CN}(\Gamma_{\mathcal{M}_{2mn}}) = \begin{cases}
	0, &\text{ if } m=3, z=1\\
	\frac{2 (m-2) (m-1) m n^2 (m n-3)}{2 m-1}, &\text{ if }2\nmid m \\
	\frac{(m-4) (m-2) m n^2 (m n-3)}{m - 1}, &\text{ if }2 \mid m.	
\end{cases}
\]

\end{cor}
\begin{proof}
We have 
\[
Z(\mathcal{M}_{2mn})=\begin{cases}
	\langle y^2 \rangle, & \text{ if }2\nmid m\\
	\langle y^2 \rangle\cup x^{\frac{m}{2}}\langle y^2 \rangle, & \text{ if }2\mid m
\end{cases}
\] 
and so  $|Z(\mathcal{M}_{2mn})| = n$ and $2n$ according as $2\nmid m$ and  $2\mid m$.
Therefore,
\[
\frac{\mathcal{M}_{2mn}}{Z(\mathcal{M}_{2mn})}\cong \begin{cases}
 	D_{2m}, & \text{ if }2\nmid m\\
 	D_m, & \text{ if }2\mid m.
\end{cases}
\]
Hence, the result follows from  Theorem \ref{D_2m theorem}.
\end{proof}

\begin{cor}\label{D_2m-cor2}
	The CNL-spectrum, CNSL-spectrum, CNL-energy and CNSL-energy of the commuting graphs of the dihedral groups, $D_{2m}=\langle x,y:~x^{2m}=y^2=1,~yxy^{-1}=x^{-1} \rangle$, $m\geq 3$, are given by
	
\centerline{	$\cnlspec(\Gamma_{D_{2m}})=\begin{cases}
		\left\lbrace (0)^{m + 1}, ((m - 1)(m - 3))^{(m - 2)}\right\rbrace, &\text{ if } 2 \nmid m\\
		\left\lbrace (0)^{(m + 1)}, ((m - 2)(m - 4))^{(m - 3)}\right\rbrace, &\text{ if } 2 \mid m,
	\end{cases}$}

\centerline{	$\cnqspec(\Gamma_{D_{2m}}) = \begin{cases}
		\left\lbrace (0)^m, (2(m - 2)(m - 3))^1, ((m - 3)^2)^{(m - 2)}\right\rbrace, &\text{ if } 2 \nmid m\\
		\left\lbrace (0)^{m}, (2(m - 3)(m - 4))^1, ((m - 4)^2)^{(m - 3)} \right\rbrace, &\text{ if } 2 \mid m,
	\end{cases}$}
 
\centerline{	$LE_{CN}(\Gamma_{D_{2m}}) = \begin{cases}
									\frac{2 (m-3) (m-2) (m-1) (m+1)}{2 m-1}, &\text{ if } 2 \nmid m\\
									\frac{(m - 4) (m - 3) (m - 2) (m + 1)}{m - 1},  & \text{ if }2\mid m
								\end{cases}$}

\noindent	and $LE^+_{CN}(\Gamma_{D_{2m}})= \begin{cases}
		\frac{2 (m-3) (m-2) (m-1) m}{2 m-1},  & \text{ if }2\nmid m\\
		\frac{(m-4) (m-3) (m-2) m}{m-1},  & \text{ if }2\mid m.
	\end{cases}$
\end{cor}

\begin{proof}
		We know that $|Z(D_{2m})| = \begin{cases}
				1, &\text{ if } 2 \nmid m\\
				2, &\text{ if } 2 \mid m
			\end{cases}$ and $\frac{D_{2m}}{Z(D_{2m})} = \begin{cases}
			D_{2m}, &\text{ if } 2 \nmid m\\
			D_{m}, &\text{ if } 2 \mid m.
		\end{cases}$
	Therefore, by using Theorem \ref{D_2m theorem} we get the required result.
\end{proof}

\begin{cor}\label{D_2m-cor4}
	The CNL-spectrum, CNSL-spectrum, CNL-energy and CNSL-energy of the commuting graphs of the groups $U_{6n}=\langle x,y : x^{2n}=y^3=1,~x^{-1}yx=y^{-1}  \rangle$ are given by
	
\centerline{	$\cnlspec(\Gamma_{U_{6n}})= \left\lbrace 0^4, (2n(2n - 2))^{2n - 1}, (n(n - 2))^{3(n - 1)} \right\rbrace,$}

\noindent	$\cnqspec(\Gamma_{U_{6n}})$
	
\quad	$= \left\lbrace ((2n - 1)(2n - 2))^1, ((2n - 2)^2)^{2n - 1}, (2(n - 1)(n - 2))^3, ((n - 2)^2)^{3(n - 1)} \right\rbrace,$

\centerline{	$LE_{CN}(\Gamma_{U_{6n}})= \frac{2}{5} (n-1) (2 n-1) (9 n+10)$ }
	
\noindent	and 
	$LE^+_{CN}(\Gamma_{U_{6n}}) = \begin{cases}
		\frac{6}{5} (n-1) (n+10) = 0, &\text{ for } n=1\\
		\frac{36}{5} (n-1) n^2, &\text{ for } n \geq 2. 							
	\end{cases}$

	
\end{cor}
\begin{proof}
	We have $Z(U_{6n})=\langle x^2 \rangle$ and $\frac{U_{6n}}{Z(U_{6n})}\cong D_6$ with $|Z(U_{6n})|= n$. Therefore we get the required result by putting $m = 3$ and $z = n$ in Theorem \ref{D_2m theorem}.
\end{proof}

\begin{cor}\label{D_2m-cor3}
	The CNL-spectrum, CNSL-spectrum, CNL-energy and CNSL-energy of the commuting graphs of the dicyclic groups,
	 $Q_{4n}=\langle x, y:~ y^{2n}=1,~x^2=y^n,~xyx^{-1}=y^{-1} \rangle$, $n\geq 2$, are given by
	 
\centerline{	 $\cnlspec(\Gamma_{Q_{4n}})=\left\lbrace (0)^{2n+1}, ((2n-2)(2n-4))^{2n-3}\right\rbrace$,}

\centerline{	 $\cnqspec(\Gamma_{Q_{4n}})=\left\lbrace (0)^{2n}, (2(2n-3)(2n-4))^1, ((2n-4)^2)^{2n-3}\right\rbrace$,}

\centerline{	 $LE_{CN}(\Gamma_{Q_{4n}})= \frac{4 (n-2) (n-1) (2 n-3) (2 n+1)}{2 n-1}$ and 
	 $LE^+_{CN}(\Gamma_{Q_{4n}})= \frac{8 (n-2) (n-1) n (2 n-3)}{2 n-1}.$}
\end{cor}
\begin{proof}
	We have $Z(Q_{4n})=\{1, x^n\}$  and $\frac{Q_{4n}}{Z(Q_{4n})}\cong D_{2n}$. Therefore, by using Theorem \ref{D_2m theorem}, we get the required result.
\end{proof}

We conclude this section with the following two results for finite non-abelian AC-group in general.
	\begin{thm}\label{AC-theorem}
	Let $G$ be a finite non-abelian AC-group and $X_1, X_2, \ldots, X_n$ be the distinct centralizers of non-central elements of $G$. Then the CNL-spectrum, CNSL-spectrum, CNL-energy and CNSL-energy of $\Gamma_G$ are given by \\
	$\cnlspec(\Gamma_G)=\{(0)^n, ((|X_i| - |Z(G)|)(|X_i| - |Z(G)| - 2))^{(|X_i| - |Z(G)| - 1)}$, where $1\leq i\leq n \},$  \\
	$\cnqspec(\Gamma_G)=\{(2(|X_i| - |Z(G)| - 1)(|X_i| - |Z(G)| - 2))^1, ((|X_i| - |Z(G)| - 2)^2)^{(|X_i| - |Z(G)| - 1)}$, where $1\leq i\leq n\},$\\
	$LE_{CN}(\Gamma_G) = n L_0 + \sum\limits_{i = 1}^n (|X_i| - |Z(G)| - 1) L_{X_i},$ where $L_0 = | 0 - \Delta_{\Gamma_G} |$ and $L_{X_i} = | (|X_i| - |Z(G)|)(|X_i| - |Z(G)| - 2) - \Delta_{\Gamma_G} |$ and \\
	$LE^+_{CN}(\Gamma_G) = \sum\limits_{i = 1}^n B_{X_i} + \sum\limits_{i = 1}^n (|X_i| - |Z(G)| - 1)  B'_{X_i}$, where $B_{X_i} = | 2(|X_i| - |Z(G)| - 1)(|X_i| - |Z(G)| - 2) - \Delta_{\Gamma_G} |$ and $B'_{X_i} = | (|X_i| - |Z(G)| - 2)^2 - \Delta_{\Gamma_G} |$. 
\end{thm}
\begin{proof}
	By \cite[Lemma 2.1]{JR},  we have
	$ \Gamma_G=\sqcup_{i=1}^n K_{|X_i| - |Z(G)|}$. 
	Here $|V(\Gamma_G)| = \sum\limits_{i = 1}^n |X_i| - |Z(G)|$ and $\tr(\cnrs(\Gamma_G)) = \sum\limits_{i = 1}^n (|X_i| - |Z(G)| - 1)(|X_i| - |Z(G)|)(|X_i| - |Z(G)| - 2)$. Therefore, $\Delta_{\Gamma_G} =\frac{\sum\limits_{i = 1}^n (|X_i| - |Z(G)| - 1)(|X_i| - |Z(G)|)(|X_i| - |Z(G)| - 2)}{\sum\limits_{i = 1}^n |X_i| - |Z(G)|}$. Hence, the result follows from the Theorem \ref{CN-LE-CNSL-LE^+_Kn}.				
\end{proof}

\begin{thm}\label{AC-cartesian-theorem}
	Let $G$ be a finite non-abelian AC-group and $X_1, X_2, \ldots, X_n$ be the distinct centralizers of non-central elements of $G$. If $A$ is a finite abelian group then the CNL-spectrum, CNSL-spectrum, CNL-energy and CNSL-energy of $\Gamma_{G\times A}$ are given by \\
	$\cnlspec(\Gamma_{G\times A})=\{(0)^n, (|A|(|X_i| - |Z(G)|)(|A|(|X_i| - |Z(G)|) - 2))^{(|A|(|X_i| - |Z(G)|) - 1)},\\ \text{ where } 1\leq i\leq n \}$,  \\
	$\cnqspec(\Gamma_{G\times A})=\{(2(|A|(|X_i| - |Z(G)|) - 1)(|A|(|X_i| - |Z(G)|) - 2))^1, ((|A|(|X_i| - |Z(G)|) - 2)^2)^{(|A|(|X_i| - |Z(G)|) - 1)},\text{ where } 1\leq i\leq n\}$,\\
	$LE_{CN}(\Gamma_G) = n L_0 + \sum\limits_{i = 1}^n (|A|(|X_i| - |Z(G)|) - 1) L_{X_i},$ where $L_0 = | - \Delta_{\Gamma_{G\times A}} |$ and $L_{X_i} = | |A|(|X_i| - |Z(G)|)(|A|(|X_i| - |Z(G)|) - 2) - \Delta_{\Gamma_{G\times A}} |$ and \\
	$LE^+_{CN}(\Gamma_G) = \sum\limits_{i = 1}^n S_{X_i} + \sum\limits_{i = 1}^n (|A|(|X_i| - |Z(G)|) - 1)  S'_{X_i}$, where $S_{X_i} = | 2(|A|(|X_i| - |Z(G)|) - 1)(|A|(|X_i| - |Z(G)|) - 2) - \Delta_{\Gamma_{G\times A}} |$ and $S'_{X_i} = | (|A|(|X_i| - |Z(G)|) - 2)^2 - \Delta_{\Gamma_{G\times A}} |$. 
\end{thm}
\begin{proof}
	We have $Z(G\times A)=Z(G)\times A$ and $X_1\times A, X_2\times A, \ldots, X_n\times A$ are the distinct centralizers of non-central elements of $G\times A$. Since $G$ is an AC-group, $G\times A$ is also an AC-group. Therefore,
	\[
	\Gamma_{G\times A} = \sqcup_{i = 1}^n K_{|X_i\times A| - |Z(G)\times A|} = \sqcup_{i = 1}^n K_{|A|(|X_i| - |Z(G)|)}.
	\]
	We have, $|V(\Gamma_{G\times A})| = \sum\limits_{i = 1}^n |A|(|X_i| - |Z(G)|)$ and $\tr(\cnrs(\Gamma_{G\times A})) = \sum\limits_{i = 1}^n (|A|(|X_i| - |Z(G)|) - 1)(|A|(|X_i| - |Z(G)|)(|A|(|X_i| - |Z(G)|) - 2))$. Therefore \\ $\Delta_{\Gamma_{G\times A}} =  \frac{\sum\limits_{i = 1}^n (|A|(|X_i| - |Z(G)|) - 1)(|A|(|X_i| - |Z(G)|)(|A|(|X_i| - |Z(G)|) - 2))}{\sum\limits_{i = 1}^n |A|(|X_i| - |Z(G)|)}$. Hence the result follows from the Theorem \ref{CN-LE-CNSL-LE^+_Kn}.
\end{proof}

\section{CNL (CNSL)-integral and CNL (CNSL)-hyperenergetic graphs}
We begin this section with the following theorem which follows from the results obtained in  Section 2.


\begin{thm} \label{consequence-1}
	Let $G$ be a finite non-abelian group.
\begin{enumerate}
\item If $G$ is isomorphic to $QD_{2^n}$ (where $n \geq 4$),  $PSL(2, 2^k)$ (where $k\geq 2$), $GL(2, q)$ (where $q >2$ is a prime power),  $A(n, \nu)$ and $A(n, p)$ then  $\Gamma_G$ is CNL (CNSL)-integral. 

\item If $\frac{G}{Z(G)}$ is isomorphic to $Sz(2)$, $D_{2m}$ and $\mathbb{Z}_p \times \mathbb{Z}_p$, then  $\Gamma_G$ is CNL (CNSL)-integral.

\item If $G$ is a $4$, $5$-centralizer finite group or a $(p+2)$-centralizer finite $p$-group then  $\Gamma_G$ is CNL (CNSL)-integral. 

\item If $G$ is isomorphic to $\mathcal{M}_{2mn}$,  $D_{2m}$, $U_{6n}$ and $Q_{4n}$ then  $\Gamma_G$ is CNL (CNSL)-integral. 

\item If $G$ is an AC-group then $\Gamma_G$ is CNL (CNSL)-integral.
\end{enumerate}	  
\end{thm}
\begin{thm}
The commuting graph of the quasihedral group $QD_{2^n}$ ($n \geq 4$) is  
\begin{enumerate}
	\item CNL (CNSL)-hyperenergetic  if $n \geq 5$.
	\item not CNL (CNSL)-hyperenergetic  if $n = 4$.
\end{enumerate}	
\end{thm}
\begin{proof}
	We have $|V(\Gamma_{QD_{2^n}})|=2^n-2$. Using  \eqref{LEcn-Kn} and Proposition \ref{Quasihedral-theorem}, we get
	\begin{align*}
			LE_{CN}(K_{|V(\Gamma_{QD_{2^n}})|} &- LE_{CN}(\Gamma_{QD_{2^n}})\\
			&= 2 \left(2^n-2-1\right) \left(2^n-2-2\right)-\frac{\left(2^n-8\right) \left(2^n-6\right) \left(2^n-4\right) \left(2^n+2\right)}{8 \left(2^n-2\right)}\\
			&= -\frac{2^{(-3 + n)} (-4 + 2^n) (100 - 7\times 2^{2 + n} + 4^n)}{-2 + 2^n}:=f_1(n).
	\end{align*}
	Note that $100 - 7\times 2^{2 + n} + 4^n=100 + 2^n(2^n-7\times 2^2)$. For $n\geq 5$, $2^n-7\times 2^2>0$. So $f_1(n) < 0$ and hence $\Gamma_{QD_{2^n}}$ is CNL-hyperenergetic. For $n=4$, $f_1(n)=\frac{1104}{7} > 0$. Hence $\Gamma_{QD_{16}}$ is not CNL-hyperenergetic. Also,
	\begin{align*} 
		LE^+_{CN}(K_{|V(\Gamma_{QD_{2^n}})|} &- LE_{CN}(\Gamma_{QD_{2^n}})\\
		&= 2 \left(2^n-3\right) \left(2^n-4\right)-\frac{2^{n-3} \left(2^n-8\right) \left(2^n-6\right) \left(2^n-4\right)}{2^n-2}\\
		&= -\frac{(-4 + 2^n)^2 (24 - 13\times 2^{1 + n} + 4^n)}{8 (-2 + 2^n))}:=f_2(n).
	\end{align*} 
	Note that $24 - 13\times 2^{(1 + n)} + 4^n = 24 - 13\times 2\times 2^n + 2^n\times 2^n= 24 + 2^n(2^n-26)>0$, for all $n\geq 5$. So $f_2(n)<0$ and hence $\Gamma_{QD_{2^n}}$ is CNSL-hyperenergetic.
	For $n = 4$,  $f_2(n)=\frac{1224}{7} > 0$. Hence $\Gamma_{QD_{2^n}}$ is not CNSL-hyperenergetic.	
\end{proof}

\begin{prop}
	The commuting graph of the projective special linear group $PSL(2, 2^k)$, $k\geq 2$, is not  CNL (CNSL)-hyperenergetic.
\end{prop}
\begin{proof}
	We have $|V(\Gamma_{PSL(2,2^k)})|=-2^k + 2^{3 k} - 1$. Using  \eqref{LEcn-Kn} and Proposition \ref{PSL-theorem}, we get
	\begin{align*}
		LE_{CN}&(K_{|V(\Gamma_{PSL(2,2^k)})|})  - LE_{CN}(\Gamma_{PSL(2,2^k)})\\ 
		&= 2 \left(8^k - 2^k - 2\right) \left(8^k - 2^k - 3\right)\\
		& - \left(15 + 2^{k + 1} + 3\times 2^{2 k + 1} - 11\times 8^k + 3\times 16^k + \frac{3 \left(-5\times 2^k + 4^k + 1\right)}{8^k - 2^k - 1}\right)\\
		&= \frac{2^k \left(-7\times 2^k+7\times 2^{3 k+1}+2^{8 k+1}+3\times 16^k-32^k-9\times 64^k+10\right)}{8^k - 2^k - 1}\\
		&= \frac{2^k \left( 7\times 2^{3 k + 1} -7\times 2^k + 2 \times 2^{8k}- 2^{5k} - 9 \times 2^{6k} + 3 \times 2^{4k} \right)}{8^k - 2^k - 1} > 0,
	\end{align*}
since $8^k - 2^k - 1 > 0$ , $ 7\times 2^{3 k + 1} -7\times 2^k > 0$ and $2 \times 2^{8k}- 2^{5k} - 9 \times 2^{6k} = 2^{8k}- 2^{5k} + 2^{8k} - 9 \times 2^{6k} > 0$ for all $k\geq 2$. Therefore, $\Gamma_{PSL(2, 2^k)}$ is not CNL-hyperenergetic.

For $k\geq 3$, we have
\begin{align*}
	LE^+_{CN}&(K_{|V(\Gamma_{PSL(2,2^k)})|})  - LE^+_{CN}(\Gamma_{PSL(2,2^k)})\\ 
	&= 2 \left(8^k - 2^k - 2\right) \left(8^k - 2^k - 3\right)\\
	& - \left(\frac{7\times 2^k + 17\times 4^k - 7}{8^k - 2^k - 1} - 13\times 8^k + 3\times 16^k + 11\times 2^{k + 1} + 4^{k + 1} + 5\right)\\
	&=\frac{2^k \left(- 3\times 2^k + 3\times 2^{2 k + 1} + 2^{8 k + 1} - 8^{k + 1} + 5\times 16^k + 32^k - 9\times 64^k - 2\right)}{- 2^k + 8^k - 1}.
\end{align*}
Note that $(- 3\times 2^k + 3\times 2^{2 k + 1} + 2^{8 k + 1} - 8^{k + 1} + 5\times 16^k + 32^k - 9\times 64^k - 2) = (6\times 2^{2 k} - 3\times 2^k) + (5\times 16^k - 8\times 8^k) + (32^k - 2) + (2\times 4^k\times 64^k - 9\times 64^k) > 0$, as $k \geq 3$. Thus,   $LE^+_{CN}(K_{|V(\Gamma_{PSL(2,2^k)})|}) - LE^+_{CN}(\Gamma_{PSL(2,2^k)}) > 0$. Therefore, $\Gamma_{PSL(2, 2^k)}$ is not CNSL-hyperenergetic.

 Also for $k=2$, $LE^+_{CN}(K_{|V(\Gamma_{PSL(2,2^k)})|}) - LE^+_{CN}(\Gamma_{PSL(2,2^k)})= \frac{380848}{59} > 0$. Therefore $\Gamma_{PSL(2, 2^k)}$ is not CNSL-hyperenergetic. 
\end{proof}

\begin{thm}
	The commuting graph of the general linear group $GL(2, q)$, where $q=p^n>2$ and $p$ is prime, is not CNL (CNSL)-hyperenergetic.
\end{thm}
\begin{proof}
		We have $|V(\Gamma_{GL(2, q)})|=q^4 - q^3 - q^2 + 1$. Using  \eqref{LEcn-Kn} and Proposition \ref{GL-theorem}, we get
		\begin{align*}
			LE_{CN}&(K_{|V(\Gamma_{GL(2, q)})|}) - LE_{CN}(\Gamma_{GL(2, q)})\\
			&= 2(q^4 - q^3 - q^2 + 1 - 1)(q^4 - q^3 - q^2 + 1 - 2) \\
			& \qquad \qquad\qquad \qquad\qquad- \frac{(q - 2) (q - 1) q^4 (q + 1) (2 q - 3) ((q - 1) q - 1)}{q^3 - q - 1}\\
			&= \frac{q^2 (q + 1) ((q - 1) q - 1) \left(q^2 (q - 1) ( 5q - 8 + 2q^2 (q - 2)) + 2\right)}{q^3 - q - 1} > 0,
		\end{align*}
		since $q \geq 3$ and  $5q - 8, \,\, q - 1, \,\, q - 2 > 0$. Therefore, $\Gamma_{GL(2, q)}$ is not CNL-hyperenergetic.
		
		Also
		\begin{align*}
			LE_{CN}&(K_{|V(\Gamma_{GL(2, q)})|})  - LE_{CN}(\Gamma_{GL(2, q)})\\
			&= \begin{cases}
				2(q^4 - q^3 - q^2 + 1 - 1)(q^4 - q^3 - q^2 + 1 - 2) \\
				\qquad\qquad -(\frac{(q-2) q (q+1) (q (q (q (q (2 q ((q-4) q+6)-7)-7)-1)+11)+2)}{q^3-q-1}), &\text{ if }q\leq 5\\				
				2(q^4 - q^3 - q^2 + 1 - 1)(q^4 - q^3 - q^2 + 1 - 2) \\
				\qquad\qquad-\frac{q (q+1) \left((q (q (q (q (2 q-11)+20)-16)+7)+8) q^3-12 q-4\right)}{q^3-q-1}, &\text{ if }q\geq 6
			\end{cases}\\
			&= \begin{cases}
			\frac{2 q^{11} - 6 q^{10} + 6 q^9 - 10 q^8 + 9 q^7 + 24 q^6 - 22 q^5 - 28 q^4 + 3 q^3 + 22 q^2 + 4 q}{q^3 - q - 1} \, > 0, &\text{ if }q\leq 5
			\vspace{.2cm }\\
			\frac{2 q^{11} - 6 q^{10} + 5 q^9 - 3 q^8 + 2 q^7 + 9 q^6 - 17 q^5 - 10 q^4 + 8 q^3 + 14 q^2 + 4 q}{q^3 - q - 1} \quad > 0, &\text{ if }q\geq 6.
		\end{cases}
		\end{align*} 
	For $3 \leq q \leq 5$, $2 q^{11} - 6 q^{10} + 6 q^9 - 10 q^8 + 9 q^7 + 24 q^6 - 22 q^5 - 28 q^4 + 3 q^3 + 22 q^2 + 4 q = ((2q) q^{10} - 6 q^{10}) + ((6q) q^8 - 10 q^8) + (24 q^6 - 22 q^5) + ((9q^3) q^4  - 28 q^4) + 3 q^3 + 22 q^2 + 4 q > 0$. Again, for $q \geq 6$, $2 q^{11} - 6 q^{10} + 5 q^9 - 3 q^8 + 2 q^7 + 9 q^6 - 17 q^5 - 10 q^4 + 8 q^3 + 14 q^2 + 4 q = ((2q) q^{10} - 6 q^{10}) + ((5q) q^8 - 3 q^8) + ((9q) q^5 - 17 q^5) + ((2q^3) q^4  - 10 q^4) + 8 q^3 + 14 q^2 + 4 q > 0$, as $q \geq 6$. Therefore $\Gamma_{GL(2, q)}$ is not CNSL-hyperenergetic. 
\end{proof}

\begin{thm}
	Let $G$ be a finite group such that $\frac{G}{Z(G)}\cong Sz(2)$. Then the commuting graph of $G$ is
\begin{enumerate}
	\item CNL-hyperenergetic if $|Z(G)| \geq 17$.
	\item not CNL-hyperenergetic if $1 \leq |Z(G)| \leq 16$.
	\item CNSL-hyperenergetic if $|Z(G)|\geq 16$.
	\item not CNSL-hyperenergetic if $1 \leq |Z(G)| \leq 15$.
\end{enumerate}	
 \end{thm}
\begin{proof}
	We have $|V(\Gamma_G)| = 3\times 5 z + 4 z = 19 z$, where $z = |Z(G)|$.  Using  \eqref{LEcn-Kn} and Theorem \ref{Suzuki-theorem}, for $z\geq 2$, we get
	\begin{align*}
		LE_{CN}(K_{|V(\Gamma_G)|}) - LE_{CN}(\Gamma_G) &= 2 (19 z-1) (19 z-2)-\frac{2}{19} (4 z-1) (z (105 z+31)-38)\\
						&= - \frac{120}{19} z (z (7 z-114)+15) := f_1(z)~~~~~\text{(say)}.
	\end{align*}
Clearly, for $z \geq 17$, $f_1(z) < 0$. Again, for $2\leq z \leq 16$, $z (7 z - 114) + 15 < 0$. Therefore,  $f_1(z) > 0$. For $z=1$, $LE_{CN}(K_{|V(\Gamma_G)|}) - LE_{CN}(\Gamma_G) = \frac{11040}{19} > 0.$ So, $\Gamma_G$ is CNL-hyperenergetic and  not CNL-hyperenergetic according as $z \geq 17$ and $1 \leq z \leq 16$.

We have 
\begin{align*}
	LE^+_{CN}(K_{|V(\Gamma_G)|})\! -\! LE^+_{CN}(\Gamma_G)\! &=\! 2 (19 z - 1) (19 z - 2) - \left( \frac{10}{19} (3 z-1) (z (28 z+45)-38) \right)\\
	&= - \frac{8}{19} (3 z (z (35 z-527)+24)+38) := f_2(z)~~~\text{(say)}.
\end{align*}
Clearly, for $z \geq 16$, $f_2(z) < 0$, as $35z - 527 > 0$. For $1 \leq z \leq 15$, $3 z (z (35 z - 527) + 24) + 38 < 0$ and so $f_2(z) > 0$. Therefore,  $\Gamma_G$ is CNSL hyperenergetic and not CNSL hyperenergetic according as  $z \geq 16$ and $1 \leq z \leq 15$.
\end{proof}
\begin{thm}\label{Theorem-con-Zp}
	Let $G$ be a finite group with $\frac{G}{Z(G)}\cong \mathbb{Z}_p\times\mathbb{Z}_p$. Then the commuting graph of $G$ is not CNL (CNSL)-hyperenergetic. 
\end{thm}
\begin{proof}
	 We have $|V(\Gamma_G)|=(p^2  - 1)z$, where $z = |Z(G)|$. By \eqref{LEcn-Kn} and Theorem \ref{Z_p-theorem}, we get
	\begin{align*}
		LE_{CN}&(K_{|V(\Gamma_G)|})  - LE_{CN}(\Gamma_G)\\
		&=\begin{cases}
			2 \left(p^2 z - z - 1\right) \left(p^2 z - z - 2\right) - 0, &\text{ if } p=2 \text{ and } z=1\\
			2 \left(p^2 z - z - 1\right) \left(p^2 z - z - 2\right)\\
			 \quad - (2 (p + 1) ((p - 1) z - 2) ((p - 1) z - 1)), &\text{ otherwise }
		\end{cases}\\
		&= \begin{cases}
			4 > 0, &\text{ if } p=2 \text{ and } z=1\\
			2p((p - 1)^2(p + 1)z^2 - 2) > 0, &\text{ otherwise},
		\end{cases}
	\end{align*}		
	since $2p((p - 1)^2(p + 1) z^2 - 2) > 0$, as $p \geq 2$ and $z \geq 1$. Hence, $\Gamma_G$ is not CNL-hyperenergetic. Again, $LE_{CN}(\Gamma_G) = LE^+_{CN}(\Gamma_G)$. Therefore, $\Gamma_G$ is also not CNSL-hyperenergetic.
\end{proof}
As a corollary of  Theorem \ref{Theorem-con-Zp} we get the following results. 
\begin{cor}
	Let $G$ be isomorphic to one of the following groups
	\begin{enumerate}
		\item $\mathbb{Z}_2\times Q_8$.
		\item $\mathbb{Z}_2\times D_8$.
		\item $\mathbb{Z}_4\rtimes \mathbb{Z}_4= \langle u,v:u^4=v^4=1,~vuv^{-1}=u^{-1}\rangle$.
		\item $M_{16}= \langle u,v:u^8=v^2=1,~vuv=u^5\rangle$.
		\item $\mathcal{SG}(16,3)= \langle u,v:u^4=v^4=1,uv=v^{-1}u^{-1},uv^{-1}=vu^{-1}\rangle$.
		\item $D_8* \mathbb{Z}_4= \langle u,v,w:u^4=v^2=w^2=1,uv=vu,uw=wu,vw=u^2wv\rangle$
	\end{enumerate}
	Then  the commuting graph of $G$ is neither CNL (CNSL)-intregal nor CNL (CNSL)-hyperenergetic. 
\end{cor}
\begin{proof}
	If $G$ is isomorphic to any of the above groups then $\frac{G}{Z(G)}\cong \mathbb{Z}_2\times \mathbb{Z}_2$. Hence the result follows from Theorem \ref{consequence-1}(b) and Theorem \ref{Theorem-con-Zp}.
\end{proof}

\begin{cor}
	The commuting graph of $G$ is not CNL (CNSL)-hyperenergetic if $G$ is a finite
	\begin{enumerate}
		\item  non-abelian group of order $p^3$, where $p$ is a prime.
		\item   $4$-centralizer group.
		\item  $5$-centralizer group.
		\item $(p+2)$-centralizer $p$-group.
		\item  non-abelian group such that the maximal size of the set of pairwise non-commuting elements of $G$ is $3$ or $4$.
	\end{enumerate}
\end{cor}
\begin{proof}
Parts	(a)-(d) are follow from Theorem \ref{Theorem-con-Zp}. Under the hypothesis in part (e), $G$ is either $4$-centralizer or $5$-centralizer finite group. Hence, the result follows from parts (b) and (c).
\end{proof}

\begin{thm}\label{Theorem-con-D2n}
Let $G$ be a finite group with $\frac{G}{Z(G)}\cong D_{2m}$, where $m\geq 2$ and $|Z(G)| = z$. Then the commuting graph of $G$ is
\begin{enumerate}
	\item  not CNL-hyperenergetic if $m = 2$ and $z \geq 1$; $m = 3$ and $z \leq 6$;  $m = 4$ and $z= 1, 2, 3$; $m = 5$ and $z= 1, 2$;  and $m = 6, 7, 8, 9, 10$ and $z = 1$. Otherwise, it is CNL-hyperenergetic.
	\item  not CNSL-hyperenergetic if $m = 2$ and $z \geq 1$;  $m = 3$ and $1 \leq z \leq 7$;  $m = 4$ and $z = 1, 2, 3$; $m = 5, 6$ and $z = 1, 2$;  and $m = 7, 8, 9, 10, 11$ and $z = 1$. Otherwise, it is CNSL-hyperenergetic.
\end{enumerate}
\end{thm}
\begin{proof}
We have $|V(\Gamma_G)|=2 m z - z$, where $z = |Z(G)|$. Let us consider the following cases.
	  
\noindent \textbf{Case 1. }  $m = 2$ and $z = 1$
	
 By \eqref{LEcn-Kn} and Theorem \ref{D_2m theorem}, we get 
 	
\[
LE_{CN}(K_{|V(\Gamma_G)|}) - LE_{CN}(\Gamma_G) = 2 (2 m z - z - 1) (2 m z - z - 2) - 0 = 4 > 0.
\]
and

\[
LE^+_{CN}(K_{|V(\Gamma_G)|} - LE^+_{CN}(\Gamma_G) = 2 (2 m z - z - 1) (2 m z - z - 2) - 0 = 4 > 0.
\]	
Therefore, $\Gamma_G$ is not CNL (CNSL)-hyperenergetic. 

\noindent \textbf{Case 2.}  $m = 2$ and $z \geq 2$
	
By \eqref{LEcn-Kn} and Theorem \ref{D_2m theorem}, we get 
\[
LE_{CN}(K_{|V(\Gamma_G)|}) - LE_{CN}(\Gamma_G) = \frac{1}{3} (36 z^2 - 24) > 0
\]
and
\[
LE^+_{CN}(K_{|V(\Gamma_G)|} - LE^+_{CN}(\Gamma_G) = \frac{4}{3} (9 z^2 - 6) > 0.
\]
Therefore, $\Gamma_G$ is not CNL (CNSL)-hyperenergetic.	

\noindent \textbf{Case 3.}  $m \geq 3$

By \eqref{LEcn-Kn} and Theorem \ref{D_2m theorem}, we get
\begin{align*}
		LE_{CN}(K_{|V(\Gamma_G)|}) - LE_{CN}(\Gamma_G) &= - \frac{2 (m - 1) m z (z (m^2 z - 2 m (z + 5) + 8) + 9)}{2 m - 1}
\end{align*}
Let $\alpha(m,z) = - 2 (m - 1) m z (z (m^2 z - 2 m (z + 5) + 8) + 9)$. Note that $(m^2 z - 2 m (z + 5) = m((m - 2)z - 10) > 0$, for all $m\geq 13$, and so $\alpha(m,z) < 0$. Also
\[
\alpha(3,z) = -12 z ((3 z - 22) z + 9) \begin{cases}
		>0,~~\text{ if } z\leq 6\\
		< 0,~~\text{ if } z\geq 7,
	\end{cases}
\]
\[
\alpha(4,z) = - 24z(9 + 8z(z - 4))  \begin{cases}
		>0,~~\text{ if } z =1, 2, 3\\
		< 0,~~\text{ if } z \geq 4,
	\end{cases}
\]
\[
\alpha(5,z) = - 120z(3 + z(5z - 14)) \begin{cases}
		>0,~~\text{ if } z =1, 2\\
		< 0,~~\text{ if } z \geq 3,
	\end{cases}
\]
\[
\alpha(6,z) = -60z( 9 + 4z(6z - 13)) \begin{cases}
		>0,~~\text{ if } z =1\\
		< 0,~~\text{ if } z \geq 2,
	\end{cases}
\]
\[ 
\alpha(7,z) = - 84z( 9 + z(35z - 62))  \begin{cases}
		>0,~~\text{ if } z =1\\
		< 0,~~\text{ if } z \geq 2,
	\end{cases} 
\]
\[ 
\alpha(8,z) = - 336 z( 3 + 8z(2z - 3)) = \begin{cases}
		>0,~~\text{ if } z =1\\
		< 0,~~\text{ if } z \geq 2,
	\end{cases}
\]
\[
\alpha(9,z) = - 144 z( 9 + z(63z - 82))  \begin{cases}
		>0,~~\text{ if } z =1\\
		< 0,~~\text{ if } z \geq 2,
	\end{cases}
\]
\[
\alpha(10,z) = - 180 z( 9 + 2z(40z - 46)) \begin{cases}
		>0,~~\text{ if } z =1\\
		< 0,~~\text{ if } z \geq 2,
	\end{cases}
\]
\[
 \alpha(11,z) = - 660 z( 3 + z(33z - 34)) < 0, \text{ for all } z \geq 1
\]
 and 

\[
\alpha(12,z) = -264 z (120 z^2 - 112 z + 9) < 0, \text{ for all } z \geq 1.
\]
Hence, $\Gamma_G$ is not CNL-hyperenergetic if  $m = 3$ and $z \leq 6$;  $m = 4$ and $z= 1, 2, 3$; $m = 5$ and $z= 1, 2$;  and $m = 6, 7, 8, 9, 10$ and $z = 1$. Otherwise it is CNL-hyperenergetic.


Now we determine whether $\Gamma_G$ is CNSL-hyperenergetic by considering the following subcases.

\noindent \textbf{Subcase 3.1.}  $m = 3$, $z = 1$ and $m = 4$, $z = 1$

By \eqref{LEcn-Kn} and Theorem \ref{D_2m theorem}, we get	
\[
LE^+_{CN}(K_{|V(\Gamma_G)|} - LE^+_{CN}(\Gamma_G) = \begin{cases}
		24 > 0, &\text{ for } m = 3, z = 1\\
		\frac{372}{7} > 0, &\text{ for } m = 4, z = 1.
	\end{cases}
\]
Therefore, $\Gamma_G$ is not CNSL-hyperenergetic.

\noindent \textbf{Subcase 3.2.}  $m = 3$ and $z \geq 2$;  $m = 4$ and $z \geq 2$; and $m \geq 5$ and $z \geq 1$

By \eqref{LEcn-Kn} and Theorem \ref{D_2m theorem}, we get
	\begin{align*}
		LE^+_{CN}&(K_{|V(\Gamma_G)|}  - LE^+_{CN}(\Gamma_G) \\
		&= \frac{2 (2 m - 1) ((2 m - 1) z - 2) ((2 m - 1) z - 1) - 2 (m -2 ) (m - 1) m z^2 (m z - 3)}{2m - 1}.
	\end{align*}
Let $\beta(m,z) = 2 (2 m - 1) ((2 m - 1) z - 2) ((2 m - 1) z - 1) - 2 (m -2 ) (m - 1) m z^2 (m z - 3)$. Then $\beta(m,z) = - 4 + 2m(4 - 9mz^2) - 6z + 24mz(1 - m) - 2z^2 + 24mz^2(1 - m) - 2m^2 z^2((m^2 - 3m + 2)z - 11m)$. Note that  $m^2 - 3m + 2 - 11m = m^2 - 14m + 2 > 0$, for all $m \geq 14$, and so  $\beta(m,z) < 0$. Also
\[
\beta(3,z) = -2 (- 10 + 75z + z^2(18z - 143))  \begin{cases}
		> 0,~~ \text{ if } z \leq 7\\
		< 0,~~ \text{ if } z \geq 8,
	\end{cases}
\]
\[ 
\beta(4,z) = -2 (- 14 + 147z + z^2(96z - 415))  \begin{cases}
		> 0,~~ \text{ if } z \leq 3\\
		< 0,~~ \text{ if } z \geq 4,
	\end{cases}
\]
\[
\beta(5,z) = -6 (- 6 + 81z + z^2(100z - 303))  \begin{cases}
		> 0,~~ \text{ if } z = 1, 2\\
		< 0,~~ \text{ if } z \geq 3,
	\end{cases}
\]
\[
\beta(6,z) = -2 (- 22 + 363z + z^2(720z - 1691))  \begin{cases}
		> 0,~~ \text{ if } z = 1, 2\\
		< 0,~~ \text{ if } z \geq 3,
	\end{cases}
\]
\[ 
\beta(7,z) = -2 (- 26 + 507z + z^2(1470z - 2827))  \begin{cases}
		> 0,~~ \text{ if } z = 1\\
		< 0,~~ \text{ if } z \geq 2,
	\end{cases}
\]
\[
\beta(8,z) = -6 (- 10 + 225z + z^2(896z - 1461))  \begin{cases}
		> 0,~~ \text{ if } z = 1\\
		< 0,~~ \text{ if } z \geq 2,
	\end{cases}
\]
\[
\beta(9,z) = -2 (- 34 + 867z + z^2(4536z - 6425))  \begin{cases}
		> 0,~~ \text{ if } z = 1\\
		< 0,~~ \text{ if } z \geq 2,
	\end{cases}
\]
\[
\beta(10,z) = -2 (- 38 + 1083z + z^2(7200z - 9019)) \begin{cases}
		> 0,~~ \text{ if } z = 1\\
		< 0,~~ \text{ if } z \geq 2,
	\end{cases}
\] 
\[
\beta(11,z) = -6 (- 14 + 441z + z^2(3630z - 4077)) \begin{cases}
		> 0,~~ \text{ if } z = 1\\
		< 0,~~ \text{ if } z \geq 2,
	\end{cases}
\]
\[
\beta(12,z) = -2 (15840 z^3 - 16127 z^2 + 1587 z - 46) < 0, \text{ for all } z \geq 1
\]
and
\[
\beta(13,z) = -2 (22308 z^3 - 20773 z^2 + 1875 z - 50) < 0, \text{ for all } z \geq 1.
\]
Hence, $\Gamma_G$ is not CNSL-hyperenergetic if  $m = 3$ and $2 \leq z \leq 7$;  $m = 4$ and $z = 2, 3$; $m = 5, 6$ and $z = 1, 2$;  and $m = 7, 8, 9, 10, 11$ and $z = 1$. Otherwise it is CNSL-hyperenergetic. Hence, the result follows.
\end{proof}
As a consequences of  Theorem \ref{Theorem-con-D2n} we get the following results.
\begin{cor}
	Suppose that $G$ is isomorphic to the metacyclic group $\mathcal{M}_{2mn}$.
	\begin{enumerate}
		\item If $m$ is even  then 
		\begin{enumerate}
			\item $\Gamma_G$ is not CNL-hyperenergetic whenever $m=4$ and $n \geq 1$; $m=6$ and $n = 1, 2, 3$; and $m=8, 10$ and $n = 1$. Otherwise, $\Gamma_G$ is  CNL-hyperenergetic. 
			\item $\Gamma_G$ is not CNSL-hyperenergetic whenever $m =4$ and $n \geq 1$; $m=6$ and $n = 1, 2, 3$;  and $m=8, 10, 12$ and $n = 1$. Otherwise, $\Gamma_G$ is  CNSL-hyperenergetic. 
		\end{enumerate}
		\item If $m$ is odd  then 
		\begin{enumerate}
			\item $\Gamma_G$ is not CNL-hyperenergetic whenever $m =3$ and $n \leq 6$; $m = 5$ and $n = 1, 2$; and $m = 7, 9$ and $n = 1$. Otherwise, $\Gamma_G$ is  CNL-hyperenergetic. 
			\item $\Gamma_G$ is not CNSL-hyperenergetic whenever $m =3$ and $n \leq 7$; $m = 5$ and $n = 1$; and $m = 7, 9, 11$ and $n = 1$. Otherwise, $\Gamma_G$ is  CNSL-hyperenergetic. 
		\end{enumerate}
	\end{enumerate}	
\end{cor}

\begin{cor}
Suppose that $G$ is isomorphic to the dihedral group $D_{2m}$.
\begin{enumerate}
\item If $m$ is even  then 
	\begin{enumerate}
	   \item $\Gamma_G$ is not CNL-hyperenergetic whenever $4\leq m\leq 10$. Otherwise, $\Gamma_G$ is  CNL-hyperenergetic. 
	   \item $\Gamma_G$ is not CNSL-hyperenergetic whenever $4\leq m\leq 12$. Otherwise, $\Gamma_G$ is  CNSL-hyperenergetic. 
	\end{enumerate}
\item If $m$ is odd  then 
\begin{enumerate}
	\item $\Gamma_G$ is not CNL-hyperenergetic whenever $3\leq m\leq 9$. Otherwise, $\Gamma_G$ is  CNL-hyperenergetic. 
	\item $\Gamma_G$ is not CNSL-hyperenergetic whenever $3\leq m\leq 11$. Otherwise, $\Gamma_G$ is  CNSL-hyperenergetic. 
\end{enumerate}
\end{enumerate}	
\end{cor}	
	
\begin{cor}
Let $G$ be a finite non-abelian group.  
\begin{enumerate}
		

		\item If $G$ is isomorphic to the generalized quartanion group of order $4n$, $Q_{4n}$ then
		\begin{enumerate}
			\item $\Gamma_G$ is not CNL-hyperenergetic when $n = 2, 3, 4, 5$. Otherwise, $\Gamma_G$ is CNL-hyperenergetic.  
			\item $\Gamma_G$ is not CNSL-hyperenergetic when $n = 2, 3, 4, 5, 6$. Otherwise, $\Gamma_G$ is CNSL-hyperenergetic.
		\end{enumerate}	
		\item If $G$ is isomorphic to $U_{6n}$ then 
		\begin{enumerate}
			\item $\Gamma_G$ is not CNL-hyperenergetic when $n \leq 6$. Otherwise, $\Gamma_G$ is CNL-hyperenergetic.  
			\item $\Gamma_G$ is not CNSL-hyperenergetic when $n \leq 7$. Otherwise, $\Gamma_G$ is CNSL-hyperenergetic.
		\end{enumerate}	
	\end{enumerate}
\end{cor}  

\begin{prop}
	Let $F=GF(2^n)$, $n\geq 2$ and $\nu$ be the Frobenius automorphism of $F$, that is $\nu(x)=x^2$, for all $x\in F$. The commuting graph of $A(n,\nu)$ is not CNL (CNSL)-hyperenergetic.
\end{prop}
\begin{proof}
	We have $V(\Gamma_{A(n,\nu)})=2^{2 n} - 2^n$. Therefore, by using \eqref{LEcn-Kn} and Proposition \ref{FA-proposition}, we get
	\begin{align*}
			LE_{CN}(K_{|V(\Gamma_{A(n,\nu)})|}) & - LE_{CN}(\Gamma_{A(n,\nu)})\\
			&= 2 \left( - 2^n + 2^{2 n} - 1\right) \left( - 2^n + 2^{2 n} - 2\right) - ( 2 \left(2^n - 2\right) \left(2^n - 1\right)^2 )\\
			&= 2 \left(2^n - 2\right) \left( 8^n - 4^n - 2 \right) > 0,
	\end{align*}
	since $2^n - 2 \geq 2$ and $8^n - 4^n - 2 \geq 46$, as $n \geq 2$. So, $\Gamma_{A(n,\nu)}$ is not CNL hyperenergetic. Again,  $LE_{CN}(\Gamma_{A(n,\nu)}) = LE^+_{CN}(\Gamma_{A(n,\nu)})$.  Therefore, $\Gamma_{A(n,\nu)}$ is not CNSL-hyperenergetic.
\end{proof}

\begin{prop}
	Let $F=GF(p^n)$, where $p$ is a prime. Then the commuting graph of $A(n,p)$ is not CNL (CNSL)-hyperenergetic.
\end{prop}
\begin{proof}
	We have $V(\Gamma_{A(n,p)})=p^{3 n} - p^n$. Therefore, by using \eqref{LEcn-Kn} and Proposition \ref{GF-proposition}, we get
	\begin{align*}
			LE_{CN}&(K_{|V(\Gamma_{A(n,p)})|}) - LE_{CN}(\Gamma_{A(n,p)})\\
			&= 2 \left(p^{3 n}-p^n-1\right) \left(p^{3 n}-p^n-2\right)-(2 \left(p^n+1\right)^2 \left(-3 p^{2 n}+p^{3 n}+p^n+2\right))\\
			&= 2 p^n(p^{2n}(p^n(p^{2n} - p^{n} - 1) +1) - 2)>0 ~~~~\text{ as }p\geq 2.
	\end{align*}
So, $\Gamma_{A(n,p)}$ is not CNL-hyperenergetic. Again,  $LE_{CN}(\Gamma_{A(n,p)}) = LE^+_{CN}(\Gamma_{A(n,p)})$. Therefore,  $\Gamma_{A(n,p)}$ is not CNSL-hyperenergetic.
\end{proof}


\noindent \textbf{Concluding Remarks:}
We observed that the commuting graphs of the groups discussed above are CNL (CNSL)-integral and this leads us to the following question.

\noindent \textbf{Question.} Which finite non-abelian groups are CNL (CNSL)-inetgral?

It is also observed that the commuting graphs of some AC-groups are CNL (CNSL)-hyperenergetic but some  are not CNL (CNSL)-hyperenergetic. Therefore, one can try to find general conditions  such that the commuting graphs of finite AC-groups are CNL (CNSL)-hyperenergetic.

\vspace{1cm}
\noindent \textbf{Acknowledgement:} 
The first author would like to thank DST for the INSPIRE Fellowship (IF200226).

\end{document}